\documentclass[a4paper,12pt]{article}
\usepackage[dvips]{epsfig}
\usepackage{amsmath,amssymb,amsbsy, amstext,amscd,amsfonts}
\usepackage{color,enumerate,euscript,graphicx,hyperref,srcltx,tikz,indentfirst}
\usepackage{tikz}
\usepackage{xcolor}
\usetikzlibrary{shapes.geometric}
\usetikzlibrary{ positioning,  shapes.geometric}
\input amssym.def
\input amssym.tex

\newtheorem{lemma}{Lemma}[section]%
\newtheorem{theorem}[lemma]{Theorem}%
\newtheorem{proposition}[lemma]{Proposition}%
\newtheorem{hypothesis}[lemma]{Hypothesis}%
\newtheorem{definition}[lemma]{Definition}%

\def\lg{\langle} \def\rg{\rangle} \def\char{\,\hbox{\rm char}\,}
\def\PSp{\hbox{\rm PSp}} \def\PSL{\hbox{\rm PSL}} \def\PG{\hbox{\rm PG}}
\def\Aut{\hbox{\rm Aut\,}} \def\soc{\hbox{\rm soc}}
 \def\cal{\mathcal} \def\B{{\cal B}}
\def\ZZ{\mathbb{Z}}\def\FF{\mathbb{F}}\def\f{\noindent}
\def\a{\alpha}\def\b{\beta} \def\di{\bigm|}  
\def\s{\sigma}\def\o{{\rm o}}\def\AGL{\hbox{\rm AGL}}\def\M{\mathrm{M}}\def\N{\mathrm{N}} \def\C{\mathrm{C}}
\def\A{\mathrm{A}}\def\K{\mathrm{K}} \def\t{\tau} \def\th{\theta}
\def\pf{\noindent{\it Proof.} } \def\Syl{\hbox{\rm Syl}}
\def\qed{\hfill $\Box$} \def\demo{\pf}\def\SL{\hbox{\rm SL}}
 
\def\GL{\mathrm{GL}} \def\AGL{\mathrm{AGL}}

\def\PSU{\mathrm{PSU}} \def\S{\mathrm{S}}\def\O{\Omega}
\def\D{\mathrm{D}}
\def\Cos{\mathrm{Cos}}
\def\rtimes{:}\def\P{\mathrm{P}}
\begin{document}
\begin{center}
{\bf\large Hamilton Cycles In Vertex-Transitive Graphs of Order $10p$}
\footnote{Corresponding author: Hao Yu\\
This work was supported in part by the National Natural Science Foundation of China(12301446) and Special Foundation for Guangxi Ba Gui Scholars and NSF of Guangxi (2021GXNSFAA220116).}
\end{center}

%\bigskip

\begin{center}
Huye Chen, Jingjian Li and  Hao Yu\\
\medskip
{\it {\small School of Mathematics and Information Science,\\
Center for Applied Mathematics of Guangxi, \\
Guangxi University,
Nanning 530004, P. R. China.}}
\end{center}

\renewcommand{\thefootnote}{\empty}%{footnote}}


\begin{abstract}
After long-term efforts, the Hamilton path (cycle) problem for connected vertex-transitive graphs of order $pq$ (where $p$ and $q$ are primes) was finally resolved in 2021, see \cite{DuH-pq}. Fifteen years ago, mathematicians began addressing this problem for graphs of order $2pq$. Among these studies, it was proved in 2012 (see \cite{KMZ12}) that every connected vertex-transitive graph of order $10p$ (where $p \neq 7$ is a prime) contains a Hamilton path, with the exception of a family of graphs that was recently confirmed in \cite{DLY-25}. In this paper, we achieve a further result: every connected vertex-transitive graph of order $10p$ (where $p$ is a prime) contains a Hamilton cycle, except for the truncation of the Petersen graph.  
  
\end{abstract}

\section{Introduction}\label{Introduction}
A {\it Hamilton path} (resp. {\it cycle}) is a simple path (resp. cycle) that visits every vertex in the graph exactly once. A simple graph that contains a Hamilton cycle is called a {\it Hamiltonian graph.}
In 1969, Lov\'asz \cite{L70} posed the question of whether there exist connected, vertex-transitive graphs without Hamilton paths. Later, in 1981, Alspach \cite{A81} further asked  whether infinitely many connected, vertex-transitive graphs lacking Hamilton cycles exist.
These two questions bridge two seemingly unrelated concepts together: the traversability and symmetry of graphs.

Until now, there are only four connected vertex-transitive graphs (having at least three vertices)  which contain no  Hamilton cycle:
the Petersen graph, the Coxeter graph and the truncation of the two graphs.

It has been shown that every connected  vertex-transitive graph (simply, say CVTG) of prime order $p$ contains a Hamilton cycle, see \cite{A81}. As for CVTGs of order $n=pq$ (where $p$ and $q$ are primes), partial results were established by
 Alspach \cite{A79,A81} for $n\in\{2p,3p\}$, and later extended to $n=5p$ in \cite{MP82}. It was finally shown by Du, Kutnar and  Maru\v si\v c that every CVTG of such order contains a Hamiton cycle, except for the Petersen graph, see  \cite{DKD1,DuH-pq}.

A natural extension is to consider CVTGs of order $n=2pq$, where $p$ and $q$ are primes.
In \cite{DM87}, Maru\v si\v c considered the Hamilton cycle problem for CVTGs of order $n=2p^2$; 
In \cite{KM08}, Kutnar and  Maru\v si\v c considered the case $n=4p$; 
In \cite{DZ2025}, the existence of Hamilton cycle for $n=6p$ has been proved at last; For the case $n=10p$, papers \cite{KMZ12} and \cite{DLY-25} proved the existence of the Hamilton path for the CVTGs of such order.
This paper focuses on the Hamilton cycle of CVTGs of order $10p$, where $p$ is a prime, culminating in the following main theorem.

\begin{theorem}\label{main}
Every connected vertex-transitive graph of order $10p$, where $p$ is a prime, contains a Hamilton cycle,  except for the truncation of the Petersen graph.
\end{theorem}

In our setting,  the quotient  graph induced by a semiregular automorphism of order $p$ may be the Petersen graph, which admits no   Hamilton cycles. This introduces additional challenges to our proof.
 To overcome these difficulties, we develop deeper methods that combine graph-theoretic and group-theoretic techniques. These approaches may also be naturally extended to other Hamilton cycle problems in graph theory.

Except for the results introduced as above, there are other results related to the Hamilton path (cycle) problem of CVGTs \cite{C98,DTY,KS09,DM87,DM85,MP83,MP82,Z15}, including a survey \cite{KM09}.

This paper is  organized as follows.
Following this introduction, Section $2$ presents the necessary notations, basic definitions, preliminary results, and the proof strategy for Theorem~\ref{main}. As a preparation, existences of Hamilton cycles of four families of graphs will be proved in Section $3$. In particular, the proof of Theorem~\ref{main} can be concluded into the non-quasiprimitive cases, that is $\Aut(X)$ contains a nontrivial normal subgroup inducing blocks of length $r$, and
lastly, cases when $r\in \{ 2, 5, p\}$ and $r\in \{ 10, 2p, 5p\}$  will be dealt with  in Sections $4$ and $5$, separately.

\section{Preliminaries}

Throughout this paper, all graphs are finite, undirected, and simple unless explicitly indicated. By $p$ we always denote a prime.
Given a graph $X$, by $V(X),~E(X)$ and $\Aut(X)$ we denote the vertex set, the edge set and the automorphism group of $X$, respectively. For any $\a\in V(X)$, let $\N(\a)$ denote the neighborhood of $\a$.
For disjoint subsets $U$ and $W$ of $V(X)$, $X\lg U \rg$ is the subgraph of $X$ induced by $U$, while $X[U,W]$ denotes the induced bipartite subgraph with bipartition $U$ and $W$.
In the case when $X\lg U \rg$ and $X[U,W]$ are regular, $d(U)$ and $d(U,W)$ denote the valency of $X\lg U\rg$ and $X[U,W]$, respectively.

Let $G$ be a transitive permutation group on a set $V$. For any subset $A\subseteq V$,
$G_A$ (resp. $G_{(A)}$) denotes the set-wise (resp. point-wise) stabilizer of $A$ in $G$.
A {\it block} (or {\it $|B|$-block}) of $G$ is a nonempty subset $B\subseteq V$ such that for every $g\in G$, either $B^g=B$ or $B^g\cap B=\emptyset$.  A block is {\it trivial} if $|B|=1$ or $B=V$, and {\it non-trivial} otherwise. Then $G$ is {\it primitive} on $V$ if it has only trivial blocks, and is {\it imprimitive} otherwise. For a block $B\subseteq V$, the collection $\mathcal{B}=\{B^g\mid g\in G\}$ forms a {\it block system} of $G$ on $V$. Obviously, if $N$ is a normal subgroup of $G$, then every $N$-orbit is a block of $G$ and the collection of $N$-orbits forms a block system, called the {\it $N$-block system}. If all the nontrivial normal subgroup of $G$ is transitive on $V$, then $G$ is said to be {\it quasiprimitive} on $V$, and is called to be {\it non-quasiprimitive} on $V$ otherwise.
A graph $X$ is called a {\it primitive graph} if every transitive subgroup of $\Aut(X)$ is primitive on $V(X)$, and $X$ is called an {\it imprimitive graph} otherwise.
For an imprimitive graph $X$, we further define $X$ to be a {\it quasiprimitive graph} if every transitive subgroup $G\leq \Aut(X)$ acts quasiprimitively on $V(X)$, and $X$ is a {\it non-quasiprimitive graph} otherwise.
In what follows, we recall some definitions and known facts.

\vskip 3mm
\f {\bf 2.1  Coset graphs and bi-coset graphs}

\vskip 3mm

Let $X$ be a $G$-vertex-transitive graph. Then $X$ is isomorphic to a coset graph of $G$. If $X$ is bipartite with its two parts corresponding to $G$-orbits, then $X$ is isomorphic to a bi-coset graph of the group $G$.
We now formally define coset graphs and bi-coset graphs in the following.

\begin{definition}\label{coset-graph}
Let $G$ be a group with nontrivial subgroups $L$, $R$ and $H$. Let $S$ be a subset of $G$, and $D$ a union of some double cosets of $L$ and $R$ in $G$.
\begin{enumerate}
  \item[{\rm{(1)}}] The coset graph $\Cos(G,H,\cup_{g\in S}HgH)$ is a directed graph with vertex set $[G:H]$ $($the right cosets of $H$ in $G$$)$, and arc set $\{(Hx,Hy)\mid yx^{-1}\in \cup_{g\in S}HgH\}$.

  \item[{\rm{(2)}}] The bi-coset graph $\rm{B}(G,L,R,D)$ is a bipartite graph with vertex set $[G:L]\cup [G:R]$, and edge set $\{\{Lg,Rdg\}\mid g\in G, d\in D\}$.
\end{enumerate}
\end{definition}

\begin{proposition}\label{coset-graph-pro}
Using the notations from Definition \ref{coset-graph}, we have the following properties$:$
\begin{enumerate}
  \item[{\rm{(1)}}] The coset graph $\Cos(G,H,\cup_{g\in S}HgH)$ is undirected if and only if $\cup_{g\in S}Hg^{-1}H=\cup_{g\in S}HgH.$
  \item[{\rm{(2)}}] The coset graph $\Cos(G,H,\cup_{g\in S}HgH)$ is connected if and only if $\lg H,S \rg=G.$
  \item[{\rm{(3)}}] The bi-coset graph $\rm{B}(G,L,R,D)$ is connected if and only if $\lg D^{-1}D \rg=G.$
\end{enumerate}
\end{proposition}

\vskip 3mm
\f {\bf 2.2  Semiregular automorphisms and quotient (multi)graphs}

\vskip 3mm
Let $m\geqslant 1$ and $n\geqslant 2$ be integers. An automorphism $\rho$ of a graph $X$ is called $(m,n)$-{\em semiregular} (in short, {\em semiregular}) if, as a permutation on $V(X)$, its cycle decomposition consists of exactly $m$ cycles, each of length $n$.
Let $\mathcal{P}$ be the set of all orbits of $\lg \rho \rg $, where $\rho$ is a semiregular automorphism.
Let $X_{\mathcal{P}}$ be the {\em quotient graph} corresponding to $\mathcal{P}$, the graph  whose vertex set is $\mathcal{P}$, with $A, B \in \mathcal{P}$ adjacent if there exist adjacent vertices $a \in A$ and $b \in B$ in $X$.
Let $X_\rho$ be the {\em quotient multigraph} corresponding to $\rho$,  the multigraph  whose vertex set is $\mathcal{P}$ and in which $A,B \in \mathcal{P}$ are joined by $d(A,B)$ edges.
Note that the quotient graph $X_\mathcal{P}$ is precisely the underlying graph of $X_\rho$.

Let $X$ be a $G$-vertex-transitive graph with $G\leq\Aut(X)$, and $N\unlhd G$ a nontrivial normal subgroup. Let $\mathcal{B}$ be the set of $N$-orbits on $V(X)$. The {\em normal quotient graph} $X_{\mathcal{B}}$ of $X$ induced by $N$ is defined with vertex set $\mathcal{B}$, and two vertices $B,B'\in \mathcal{B}$ are adjacent if and only if there exist vertices $\b\in B$ and $\b'\in B'$ that are adjacent in the graph $X$.

\vskip 3mm
\f {\bf 2.3  Lifting cycle technique}
\vskip 3mm

A key tool for addressing Hamilton cycle problems is the lifting cycle technique (see \cite{A1989,KM09,DM83}).
When the quotient graph is applied relative to a semiregular automorphism of prime order and the corresponding quotient multigraph possesses two adjacent orbits linked by a double edge encompassed within a Hamilton cycle, lifts of Hamilton cycles from quotient graphs are invariably achievable. This double edge enables us to conveniently ``change direction" to procure a walk in the quotient graph that elevates into a full cycle above.

Let $X$ be a graph admitting a $(m,n)$-semiregular automorphism $\rho$.
Let $\mathcal{P} = \{S_1, S_2,\\ \cdots , S_m\}$ be the set of all orbits of $\langle\rho\rangle$, 
and $\pi : X \to X_{\mathcal{P}}$ the corresponding projection of $X$ to its quotient graph $X_{\mathcal{P}}$. For a (possibly closed) path $W = S_{i_1}S_{i_2}\ldots S_{i_k}$ in $X_{\mathcal{P}}$ we let the {\em lift} of $W$ be the set of all paths in $X$ that project to $W$. The proof of the following lemma is straightforward and is just a reformulation of \cite[Lemma~5]{MP82}.

\begin{proposition} \label{pro:4}
Let $X$ be a graph admitting
a $(m,p)$-semiregular automorphism $\rho$, where $p$ is a prime.
Let $C$ be a cycle of length $k$ in the quotient graph $X_{\mathcal{P}}$,
where $\mathcal{P}$ is the set of orbits of $\langle\rho\rangle$.
Then, the lift of $C$ either contains a cycle of length
$kp$ or consists of $p$ disjoint $k$-cycles.
In the latter case, we have $d(S,S') = 1$ for every edge $SS'$ of $C$.
\end{proposition}

\vskip 3mm
\f {\bf 2.4 Some results  on graph theory}
\vskip 3mm
Here is a particularly useful conclusion regarding Hamilton cycles that follows directly from \cite[Theorem~1]{BJ78}.

\begin{proposition} \label{pro:2}
Every connected vertex-transitive graph of order $n$ with valency at least $n/3$ contains a Hamilton cycle.
\end{proposition}

For integers $n\geqslant 2$ and $1\leq k\leq n-1$,
the generalized Petersen graph $\operatorname{GP}(n, k)$ has vertex-set $\{u_i,v_i\di i\in\ZZ_n\}$ and edge-set
$\{u_i u_{i+1}, u_i v_i, v_i v_{i+k}\di i\in\ZZ_n\}$.
By combining Theorem $3$ from \cite{A1989} with the result on page $119$ of \cite{Biggs}, we establish the following proposition.
\begin{proposition} \label{petersen}
Every vertex-transitive generalized Petersen graph either contains a Hamilton cycle or is isomorphic to the Petersen graph.
\end{proposition}
\begin{proposition}\label{10pcayley}{\rm \cite[Theorem~1.2]{M2022}}
Let $G$ be a group of order $10p$, where $p$ is a prime.
Then every connected Cayley graph on $G$ contains a Hamilton cycle.
\end{proposition}

In what follows, we review key known results and technical properties.
\begin{proposition} \label{pq-H-cycle}{\rm \cite[Theorem~1.4]{DuH-pq}}
With the exception of the Petersen graph, a connected vertex-transitive graph of order $pq$ with $p$ and $q$ being primes, contains a Hamilton cycle.
\end{proposition}

\begin{proposition} {\rm \cite[Remark~3.9]{DZ2025}}\label{quasi-pri-2} Let $X$ be a $G$-vertex-transitive graph of order $10p$, where $p\ge7$ is a prime.
Let $N$ be a normal subgroup of $G$ and $N\cong\SL(2,2^{2^s})$ induces $N$-blocks of size $5p$, where $p=1+2^{2^s}$ and $5\di (p-2)$.
Suppose that $N$ is faithful, imprimitive and quasiprimitive on both $N$-blocks. Then $X$ contains a  Hamilton cycle.
\end{proposition}

\begin{proposition} \label{pro:3}  {\rm \cite[Lemma~2.4]{KS09}}
    Let $X$ be a $G$-imprimitive graph of order $mq$, where $q$ is a prime, and let $G$ be an imprimitive subgroup of $\Aut(X)$ and $N$ be a normal subgroup of $G$ with orbits of length $q$. Then $X$ has a $(m,q)$-semiregular automorphism whose orbits coincide with the orbits of $N$.
\end{proposition}

\begin{proposition} {\rm \cite[Theorem~2.4]{A1989}}\label{deg3}
Let $X$ be a graph that admits a semiregular automorphism $\rho$ of order $t \geqslant 4$ and $\mathcal{P}$ the set of $\lg\rho\rg$-orbits on $V(X)$. Let $X_1, X_2, \ldots, X_m$ be the subgraphs of $X$ induced by the orbits of $\lg\rho\rg$,
and let each $X_i$ be connected.
If $d(X_j) \geqslant 3$ for some $j$ and there is a Hamilton path in the quotient graph $X_{\mathcal{P}}$ with $X_j$ as one of its endvertices, then $X$ contains a Hamilton cycle.
\end{proposition}

\begin{proposition} {\rm \cite[Lemma~2.7]{DZ2025}}\label{direct product}
 Let $H$ be a subgroup of $T$ such that $|T:H|=t$ and $G =T\times \lg c \rg$, where $|c|=p$ is a prime such that $(p,t)=1$. Let $D=HgH\cup Hg^{-1}H$ for some $g\in T\backslash H $ and set $Y_1=\Cos(G,H,Dc\cup Dc^{-1})\cong \Cos(G,H,Dc^k\cup Dc^{-k})$ and $Y_2 = \Cos(G,H,D\cup HcH \cup Hc\sp{-1}H)\cong \Cos(G,H,D\cup Hc^kH \cup Hc\sp{-k}H)$, where $k\neq 0$. Then we have$:$
\begin{enumerate}
\item[{\rm{(1)}}] suppose that the graph $X=\Cos(T,H ,D)$ contains a Hamilton cycle.
  Then both $Y_1$ and $Y_2$ contain a Hamilton cycle$;$
\item[{\rm{(2)}}] suppose that $T$ acts primitively on the set $[T:H]$ and $|G:H|=2rs$ with two distinct odd primes $r$ and $s$.
  Then for any subset $D_0$, where $\lg D_0\rg=G$, there exists some $D'=Hg'H\cup Hg'^{-1}H$, where $g'\in T\setminus H$, such that  the connected coset graph $\Cos(G, H, D_0)$ contains a subgraph with the form either $\Cos(G, H, D'c^k\cup D'c^{-k}))$ or $\Cos(G, H, D'\cup Hc^kH \cup Hc\sp{-k}H)$, where $k\neq 0.$
\end{enumerate}		
\end{proposition}

\vskip 3mm
\f {\bf 2.5 Some results  on group theory}
\vskip 3mm
Group theory plays a fundamental role in our work. We begin by compiling relevant group-theoretic results that will be essential for our analysis. Throughout this paper, we employ standard notation and terminology as found in \cite{Dixon} and other canonical references.
We begin by recalling fundamental results on subgroup complements. Let $G$ be a group and $H$  a subgroup of $G$. Then a {\em complement} of $H$ in $G$ is a subgroup $K$ of $G$ such that $G=HK$ and $H\cap K=1$.

\begin{proposition}{\rm\cite[Satz 1]{G1952}}\label{complement}
Let $N\leq M\leq G$ with $(|N|,|G:M|)=1$ and $N$ an abelian normal subgroup of $G$. If $N$ has a complement in $M$, then $N$ also has a complement in $G$.
\end{proposition}

\begin{proposition}{\rm\cite[3.3.1 Theorem of Schur-Zassenhaus]{KS}}\label{complement2}
Let $K$ be an abelian normal subgroup of $G$ with $(|K|,|G:K|)=1$. Then $K$ has a complement in $G$, and all the complements of $K$ are conjugate in $G$.
\end{proposition}

\begin{proposition}\label{the:main5} {\rm \cite[p. 42 Theorem~7.8]{fgt}}
Suppose that $G$ is a finite group with a normal subgroup $H$ and $P$ is a $Sylow$ $p$-subgroup of $H$ $($i.e., $P\in\Syl_p(H)$$)$. Then $G=\N_G(P)H$.
\end{proposition}

\begin{proposition} \label{the:main4} {\rm \cite[Theorem~1.49]{fsg}}
Every transitive permutation group $G$ of prime degree $p$ is one of the following cases $($up to isomorphic$)$$:$ $\ZZ_p\rtimes \ZZ_r\lesssim\AGL(1,p)$ with $r\di (p-1)$ and  
$(\soc(G),p)\in\{(A_p,p),(\PSL(2,11),11),(\M_{11},11),(\M_{23},23),(\PSL(d,q),\frac{q\sp d-1}{q-1})\}$.
\end{proposition}

The following result is quite trivial, and we omit its proof.
\begin{lemma}\label{deg10}
There are totally 5 minimal transitive groups with degree 10 up to isomorphic:
$\ZZ_{10},D_{10},\ZZ_5^2\rtimes\ZZ_{2^m},\ZZ_2^4\rtimes\ZZ_5$,
$\ZZ_5\rtimes\ZZ_4$ and $\A_5$, where $m$ is a positive integer.
\end{lemma}

The following two results can be extracted from \cite[Corollary 1.2]{LS85} and \cite[Theorem 1.2]{LL}.
\begin{proposition}\label{2p}
Let $G$ be a primitive permutation group with degree $2p$, where $p>5$ is a prime. Then $G$ is isomorphic to $\M_{22}$ with degree $22$, $\A_{2 p}$ with degree $2p$ or $\PSL(2, r^{2^i})$ with degree $r^{2^i}+1$, where $r$ is an odd prime and $i$ is a positive integer.
\end{proposition}

\begin{proposition}\label{5p}
Let $G$ be a primitive transitive permutation group of degree $5p$ for $p\geqslant 7$ prime. 
Then, up to isomorphic, $(\soc(G), p)\in\{(\A_{5p}, p)$, $(\PSL(d,q), \frac{q^d-1}{5(q-1)})$, $(\P\Omega(5,4), 17)$, $(\PSU(3,4), 13)$, $(\A_{11}, 11),(\A_7, 7),(\PSL(4,2), 7),(\PSL(5,2), 31),(\PSp(4,4), 17),(\P\Omega^+(6, 2),$ $7),(\PSL(2, 5^2),$ $ 13)$, $(\M_{11}, 11)$, $(\PSL(2, 11), 11)\}$.
Moreover, $G$ is $2$-transitive in the first four cases, and $\soc(G)$ has a metacyclic transitive subgroup isomorphic to $\ZZ_p:\ZZ_5$ if
$\soc(G)\in\{\A_{11}, \PSL(5,2), \M_{11}, \PSL(2,11)\}.$
\end{proposition}

The following result concerns quasiprimitive and imprimitive permutation groups whose degree is a product of two distinct primes, as established in \cite[p. 313]{PWX}.
\begin{proposition}\label{quasiprimitive}
Let $G$ be a quasiprimitive and imprimitive permutation group of degree $pq$, where $p\geqslant 7$ is prime and $q\in\{2,5\}$. Then $G$ is an almost simple group satisfying one of the following:

\begin{itemize}
  \item [{\rm(1)}] when $q=2:$
   \begin{itemize}
     \item[{\rm(a)}]$\soc(G)\cong\M_{11}$ with $p=11,$ or
     \item[{\rm(b)}]$\soc(G\cong\PSL(d,e)$ for some odd primes $d$ and $e$ with $(d,e-1)=1$ and $p=\frac{e^d-1}{e-1};$
   \end{itemize}
  \item [{\rm(2)}] when $q=5:$
  \begin{itemize}
    \item[{\rm(a)}] $\soc(G)\cong\SL(2,2^{2^s})$ with $p=2^{2^s}+1$ and $5|2^{2^s}-1,$ or
    \item[{\rm(b)}] $\soc(G)\cong\PSL(d,e)$ for odd primes $d$ and $e$ with $(d,e-1)=1,$ $p=\frac{e^d-1}{e-1}$ and $5\mid e-1,$ or
    \item[{\rm(c)}] $\soc(G)\cong\PSL(2,11)$ with $p=11.$
  \end{itemize}
\end{itemize}

\end{proposition}

Let $G$ be a permutation group acting on a set $\Omega$.
If $H \leq G$ is a subgroup and $A \subseteq \Omega$ is a subset invariant under the action of $H$, then we let $H^A$ denote the restriction of $H$ to $A$.

\begin{proposition} {\rm \cite[Theorem 2.10]{DMMN2007}}\label{squarefree}
Let  $G$ be a transitive permutation group of square-free degree $n$, with a complete block system $\B$ having blocks of size  $m$  formed by the orbits of a minimal normal subgroup $N$ of $G$.
Then for all blocks  $B\in\B$ precisely one of the following occurs.
\begin{enumerate}
    \item  $ N^B $  is simple and thus quasiprimitive$;$ or
    \item  $ m = 30 $  and  $ N^B \cong \A_5^2 $ , with one copy of  $ \A_5 $  acting on $6$ blocks of size $5$ and the other copy of  $ \A_5 $  acting on $5$ blocks of size $6;$ or
    \item  $ m = 105 $  and  $ N^B \cong \A_7^2 $ , with one copy of  $ \A_7 $  acting on $15$ blocks of size $7$ and the other copy of  $ \A_7 $  acting on $7$ blocks of size $15;$ or
    \item  $ m = 6090 $  and  $ N^B \cong \PSL(2, 29)^2 ,$ with one copy of  $\PSL(2, 29) $  acting on $30$ blocks of size $203$ and the other copy of  $\PSL(2, 29) $  acting on $203$ blocks of size $30.$
\end{enumerate}
\end{proposition}

\vskip 3mm
\f {\bf 2.6 Strategy of the proof of Theorem~\ref{main}}
\vskip 3mm
To prove Theorem~\ref{main}, let $X$ be a connected vertex-transitive graph of order $10p$, where $p$ is prime.
By the main results of \cite{DZ2025,KM08,KM09} and \cite{DM87}, all connected vertex-transitive graphs of order $4p$, $2p^2$ or $6p$ are Hamiltonian, except for the truncation of the Petersen graph. Consequently, we may restrict our consideration to the case where $p\ge 7$.

From \cite{KMZ12} and \cite{DLY-25} (resp. \cite{DTY}), $X$ is  Hamiltonian, provided that  $X$ is  quasiprimitive (resp. primitive).
Henceforth, we assume that $X$ is non-quasiprimitive. Reasonably, we make the following hypothesis.
\begin{hypothesis} \label{hypothesis1}
Let  $X$ be  a connected vertex-transitive  graph of order $10p$,  where   $p\geqslant7$ is a prime, and $G$  a minimal transitive subgroup of $\Aut(X)$, which contains a nontrivial normal subgroup $N$ inducing a nontrivial complete block system, saying   $\B:=\{B_0,B_1,B_2,\cdots, B_{m-1}\}$, while by $K$ we denote the kernel of $G$ on $\B$.
Let $X_\B$ be the quotient graph of order $m=\frac {10p}r$, where $r$ might be $2$, $5$, $p$, $10$, $2p$ or $5p$.
\end{hypothesis}

We remark that $G/K$ is a minimal vertex transitive subgroup of $\Aut(X_\B)$, following the minimality of $G$.

In the coming Section 3,   we shall show that  four spacial  families of graphs are hamiltonian; in Sections 4 and 5, we shall deal with cases when $r=|B_i|\in \{2, 5, p\}$  and $r\in \{ 10, 2p, 5p\}$, respectively, and in all cases, $X$ will be shown to be  hamiltonian, see  Theorem \ref{p-blocks}, \ref{q=5}, \ref{q=2},  \ref{l=10}, \ref{l=5p} and \ref{l=2p}, respectively.
Thus Theorem~\ref{main} is proved.

\section{Four families of graphs}\label{special}
To facilitate subsequent proofs, we shall show the following four families of graphs are Hamiltonian. 
Under Hypothesis \ref{hypothesis1},
\vskip 3mm
(i) firstly, in Lemmas~\ref{quasi-pri-1} and \ref{quasi-pri-2p},  we will deal with two families of graphs such that $\soc(K)\cong\PSL(d,e)$ acting on $r$-blocks, where  $r=5p$ or $2p$, respectively; and
\vskip 3mm
(ii) next,  replacing $5$ by any prime $q$, in Lemmas \ref{unfaithful1} and  \ref{unfaithful2}, we will deal with  another two  families of graphs of order $2pq$  where $q$ and $p$  are two  odd primes such that $q<p$ and either the kernel $K$ acts unfaithfully on a $r$-block, where $r\in \{ 2, p, q\}$; or  there exists a minimal normal subgroup $N_0$ of $G$ acts unfaithfully on a $r$-block, where $N_0$ induces the same blocks with $N$ and $r\in\{2p,2q,pq\}$.

 \subsection{Two families of graphs with $\soc(K)\cong \PSL(d,e)$}
\begin{lemma} \label{quasi-pri-1} Under Hypothesis \ref{hypothesis1}, suppose that $\soc(K)\cong\PSL(d,e)$ is faithful, imprimitive and quasiprimitive on both two blocks of size  $5p$, where $p=\frac{e^d-1}{e-1}$, both $d$ and $e$ are primes, $(d,e-1)=1$ and $5\di (e-1)$. Then $X$ is Hamiltonian.
\end{lemma}
\demo Suppose that $T:=\soc(K)\cong\PSL(d,e)\cong\SL(d,e)$, as $(d, e-1)=1$. Then $T\le N\le K\le \Aut(T)$
and  $V(X)$ consists of two $N$-blocks of size $5p$. Clearly, $T^g=T$ for any  element $g\in G\setminus K$.
To show $X$ is hamiltonian, it suffices to examine the action of the group $T\lg g\rg$ on the graph $X$. Moreover, without loss of generality, we may assume $N=T\cong \PSL(d,e)\cong \SL(d,e)$.

Recall that the projective points of $V(d,e)$ are all one-dimensional subspaces of $V(d,e)$, while the hyperplanes of $V(d,e)$ comprise all $(d-1)$-dimensional subspaces.
We note that $N$ admits two non-equivalent permutation representations of degree $p=\frac{e^d-1}{e-1}$: one given by its action on the projective points, and another by its action on the hyperplanes.
Set $\lg \theta \rg=\FF_{e}^*$. Take a point and hyperplane in $\PG(V)$:
 $$v=\lg (1,0,0,\cdots,0)\rg_{\FF_{e}}, \quad w=\{\lg (0,x_2,x_3,\cdots,x_{d})\rg_{\FF_{e}}\mid x_i\in\FF_{e}, 2\leq i\leq d\}.$$
Then  $$N_v=\{\left(\begin{array}{cc}
                                x & 0 \\
                                \beta^T & B \\
                              \end{array}\right)_{d\times d}\mid x\in\FF_{e}^*, \beta\in V(d-1, e), B\in\GL(d-1, e), x|B|=1 \}$$
and $$N_w=\{ \left(\begin{array}{cc}
                                x & \beta \\
                                0 & B \\
                              \end{array}\right)_{d\times d}\mid x\in\FF_{e}^*, \beta\in V(d-1, e), B\in\GL(d-1, e), x|B|=1 \}.$$

\f Set
$$N_1=\{\left(
                              \begin{array}{cc}
                                1 & 0 \\
                                \beta^T & E_{d-1} \\
                              \end{array}
                            \right)_{d\times d}\mid \beta\in V(d-1, e)\},
\quad N_2=\{\left(
                              \begin{array}{cc}
                                1 & \beta \\
                                0 & E_{d-1} \\
                              \end{array}
                            \right)_{d\times d}\mid \beta\in V(d-1, e)\},$$$$H=\{\left(
                              \begin{array}{cc}
                                \theta^{5i} & 0 \\
                                0 & B \\
                              \end{array}
                            \right)_{d\times d}\mid  B\in \GL(d-1,e), |B|=\theta^{-5i}, i\in\ZZ\} \,{\rm and}\,
t=\left(
                              \begin{array}{ccc}
                                \theta & 0&0 \\
                                0 & \theta^{-1}& 0 \\
                                0&0&E_{d-2}
                              \end{array}
                            \right).$$
A straightforward verification confirms that $N_1$, $N_2$, $H$, $N_1H$ and $N_2H$ are all groups, in particular, $N_1\cong \ZZ_e^{d-1}\cong N_2$ and $t^5\in H$.  Furthermore, we have the equalities  $N_{v}=N_1H\lg t\rg$ and $N_w=N_2H\lg t\rg$.

Define $L=N_1H$ and $R=N_2H$. Then $|N_v:L|=|N_w:R|=5$. Since $|N:N_v|=|N:N_w|=p$, employing the \cite[Lemma 4.5]{PWX}, one yields that the imprimitive and quasiprimitive action of $N$ on $5p$ points is permutation equivalent to the action of the group $N$ on the right cosets $[N:L]$ or $[N:R]$.

Let $\Gamma=\rm{B}(N, A_1, A_2, D)$ be a bi-coset graph, where   $D$ is a union of some double cosets of $A_1$ and $A_2$ in $N$  for some   $A_1,A_2\in\{ L, R\}$. Then $\Gamma$ is a subgraph of graph $X$. Since a Hamilton cycle in any generating subgraph of $X$ automatically yields a Hamilton cycle in $X$ itself, it suffices to consider the subgraph $\Gamma$.

Let $\FF_{e^d}=\FF_e(\alpha)$ be the extension of $\FF_e$ by adding a root $\alpha$ of the irreducible polynomial $f(x)=\sum_{i=0}^d a_ix^i$ over $\FF_{e}$. Consequently, every element in $\FF_{e^d}^*$ of order $e^d-1$ can be expressed in the form  $\sum_{i=0}^{d-1} s_i'\a^i$. Moreover, $\FF_{e^d}$ can be identified with the vector space $V(d, e)$, where each element $a=\sum_{i=0}^{d-1} s_i'\a^i\in\FF_{e^d}$ induces a linear transformation $\rho_a: x\mapsto ax$, for all $x\in \FF_e$. Under the basis
$\{1, \a, \cdots, \a^{d-1}\}$, the linear transformation $\rho_a$ satisfies:
 $$\rho_a((1, \a, \a^2, \cdots, \a^{d-1})^T)=A'(1, \a, \a^2, \cdots, \a^{d-1})^T,$$
where  $A':=\left(\begin{array}{cccc}s_0' & s_1' & \cdots & s_{d-1}' \\\ast & \ast  & \cdots  & \ast \\
\end{array}\right)_{d\times d}$. 
This embeds $\FF_{e^d}^*$ in $\GL(d,e)$, while its determinant-$1$ matrices form a cyclic subgroup of $\SL(d,e)$ of order $p=\frac{e^d-1}{e-1}$, generated by  $A:=\left(\begin{array}{cccc}s_0 & s_1 & \cdots & s_{d-1} \\
\ast & \ast  & \cdots  & \ast \\\end{array}\right)_{d\times d}$ for some $s_j\in\FF_e$ and $j\in\{0, 1, \cdots, d-1\}$.

Note that $L\leq N_v\leq N$ and $R\leq N_w\leq N$. Then the sets $\B_v=\{[N_v:L]g\mid g\in N\}$ and $\B_w=\{[N_w:R]g\mid g\in N\}$ form two distinct $N$-block systems for the right multiplication action of $N$ on $[N:L]$ and $[N:R]$, respectively. Let $\mathcal{C}=\mathcal{B}_1\cup \mathcal{B}_2$ with $\mathcal{B}_1,\mathcal{B}_2\in\{\mathcal{B}_v,\mathcal{B}_w\}$. Then $|\B_1|=|\B_2|=p$ and $X_\mathcal{C}$ is a $5$-blocks quotient graph of $X$ induced by the group $\lg t\rg$, where $t$ is  semiregular.

Observe that $N\cong\SL(d, e)$ has exactly two conjugacy classes of subgroups of index $p$ in $N$. Consequently, $N$ admits two nonequivalent permutation presentations of degree $p$: the transitive permutation presentation of $N$ on projective points ($[N:N_v]$) and projective hyperplanes ($[N:N_w]$). Moreover, these two transitive permutation presentations are $2$-transitive. 
Then we may assume that $\{\B_1,\B_2\}\subseteq\{\mathcal{B}_v,\mathcal{B}_w\}$, which leads to four possible isomorphism types for the bipartite subgraph $X[\mathcal{B}_1,\mathcal{B}_2]$: $(1)$ the non-incidence graph of points and hyperplanes in $PG(d-1,e)$, $(2)$ the incidence graph of points and hyperplanes in $PG(d-1,e)$,  $(3)$ $K_{p,p}-pK_{1,1}$ or $K_{p,p}$ and $(4)$ $pK_{1,1}$. We now analyze these four cases separately.

\vskip 3mm
{\bf Case 1}: $X[\mathcal{B}_1,\mathcal{B}_2]$ is the non-incidence graph of points and hyperplanes in $\PG(d-1,e)$.
\vskip 3mm

Now consider the bi-coset graph $\Gamma=\rm{B}(N, L, R, RL)$. By the  definition of $\lg A\rg$, the element $A$ is a $(5,p)$-semiregular automorphism of $\Aut(\Gamma)$. Then the quotient graph $X_{\mathcal{B'}}$ of $\Gamma$ induced by $\lg A\rg$  is a bipartite graph of order $10$. In what follows,  we assert that the quotient graph $X_{\mathcal{B'}}$ contains a Hamilton cycle, which can be lifted to that of the graph $X$.

Recall that $\lg A \rg $ is regular on both projective points and  projective hyperplanes, and the definition of $H$. Then it is easy to see that $H$ has five orbits when acting on the set $\{(x,0,0,\cdots,0)\mid x\in\FF_{e}^*\}$ and hence the vertex set of $X_{\mathcal{B'}}$ can be re-labeled as $$V(X_{\mathcal{B'}})=\{B_{k,0}\di k\in \ZZ_5\}\cup \{B_{k,1}\di k\in \ZZ_5\},$$  where 
$$B_{k,0}=\{N_1Ht^ks\mid s\in\lg A \rg\}\quad{\rm and }\quad B_{k,1}=\{N_2Ht^ks\mid s\in\lg A \rg\}.$$
Note that $\N(L):=\{Rl\mid l\in L\}$ is the neighborhood of the vertex $L$ in the bi-coset graph $\Gamma$. Then $|\N(L)|=\frac{|L|}{|L\cap R|}=\frac{|N_1H|}{|N_2H\cap N_1H|}=e^{d-1}$, i.e., the valency of $\Gamma$ is $e^{d-1}$, as
$$N_1H=\{\left(
                              \begin{array}{cc}
                                \theta^{5i} & 0 \\
                               \beta^T  & B \\
                              \end{array}
                            \right)_{d\times d}\mid \beta\in V(d-1,e), B\in \GL(d-1,e), |B|=\theta^{-5i}, i\in\ZZ\}$$
and $N_1H\cap N_2H=H.$
Then we carry out our proof of Case 1 by three steps below.

\vskip 3mm

{\it Step 1.1: Show that $\Aut(X_{\mathcal{B'}})$ contains a $(2,5)$-semiregular automorphism $\phi'$, where $\lg\phi'\rg$ induces two orbits $\{B_{k,0}\mid k\in\ZZ_5\}$ and $\{B_{k,1}\mid k\in\ZZ_5\}$ on the set $V(X_{\mathcal{B'}})$.}
\vskip 3mm

Let $\N(Lt^\ell)$ be the neighborhood  of the vertex $Lt^\ell$ in the bi-coset graph $\Gamma=\rm{B}(N, L, R, \\RL)$, where $0\leqslant \ell\leqslant 4$.
Note that $Lt^\ell\in B_{\ell, 0}$ and $\N(Lt^\ell)=RLt^\ell:=\{Rlt^\ell\mid l\in L\}$. 
Then $B_{k,1}\in \N(B_{\ell,0})$ if and only if $Rt^k\iota\in \N(Lt^\ell)$ for some $\iota\in \lg A \rg$. Note that $H^{N_\kappa}=H$, $H^t=H$ and $N_\kappa^t=N_\kappa$ for $\kappa=1, 2$. Recall that  $N_v=H:N_1:\lg t\rg$ and $N_w=H:N_2:\lg t\rg$ with $L=H:N_1$ and $R=H:N_2$. Therefore, we deduce that $\N(Lt^\ell)=Rt^\ell N_1:=\{Rt^\ell n_1\mid n_1\in N_1\}$
and $B_{k,1}\in \N(B_{\ell,0})$ if and only if $Rt^k\lg A\rg\cap Rt^\ell N_1\neq \emptyset$, i.e., $Rt^\ell n_1=Rt^k \iota$ for some $n_1\in N_1$ and $\iota\in \lg A \rg$. Equivalently, there exist $\iota'\in\lg A\rg$ and $n'\in N_1$ such that $Rt^{k-\ell}\iota'n'=R$, i.e., $t^{k-\ell}\iota'n'\in R$.
Without loss of generality, set $n'=\left(
  \begin{array}{cc}
    1 & 0 \\
    \eta^T & E_{d-1} \\
  \end{array}
\right)$ and $\iota'=\left(
  \begin{array}{cc}
    \theta^{\delta_1} & \gamma \\
    (\gamma')^T & S \\
  \end{array}
\right),$ where $S$ is a submatrix of $\iota'$ such that $\iota'\in \lg A\rg$ and $\eta=(\beta_1,\beta_2,\cdots,\beta_{d-1})$,
$\gamma=(s_1,s_2,\cdots,s_{d-1})$ and $\gamma'=(-s_{d-1}a_{d}^{-1}a_0,-s_{d-2}a_{d}^{-1}a_0,\cdots,-s_{1}a_{d}^{-1}a_0)$
 for some $\theta^{\delta_1},\beta_i,s_i\in\FF_e$   with $1\leq i\leq d-1$.
Further by the definition of $t$ and $R$, combing with a simple computation, one yields that
$$\begin{array}{lll}
    t^{k-\ell}\iota' n'\in R &\Longleftrightarrow &{\tiny \left(\begin{array}{ccccc} \theta^{k-\ell} & 0 & 0 \\
0 & \theta^{-(k-\ell)} & 0 \\
 0 & 0 & E_{d-2}  \\
 \end{array}\right)
\left(
  \begin{array}{cc}
    \theta^{\delta_1} & \gamma \\
    (\gamma')^T & S \\
  \end{array}\right)}{\tiny\left(
  \begin{array}{cc}
    1 & 0 \\
    \eta^T & E_{d-1} \\
  \end{array}
\right)\in R}\\
    &\Longleftrightarrow& {\tiny \left(
  \begin{array}{cc}
    \theta^{\delta_1} & \gamma \\
    (\gamma')^T & S \\
  \end{array}\right)\left(
  \begin{array}{cc}
    1 & 0 \\
    \eta^T & E_{d-1} \\
  \end{array}
\right)=}{\tiny \left(\begin{array}{ccccc} \theta^{\ell-k} & 0 & 0 \\
0 & \theta^{k-\ell} & 0 \\
 0 & 0 & E_{d-2}  \\
 \end{array}\right)\left(
                                \begin{array}{cc}
                                  \theta^{5i'} & \beta \\
                                  0 & B \\
                                \end{array}
                              \right)}\\
                               &&  \text{({ for some }{\tiny $\left(
                                \begin{array}{cc}
                                  \theta^{5i'} & \beta \\
                                  0 & B \\
                                \end{array}
                              \right)\in R $})} \\    
 \end{array}  $$                           
$$\begin{array}{lll}                                                 &\Longleftrightarrow&
    {\tiny \left(
\begin{array}{ccccc}
 \theta^{\delta_1} &  s_1 &s_2 & \cdots & s_{d-1} \\
-s_{d-1}a_d^{-1}a_0 & \ast  &\ast  & \cdots  & \ast \\
 -s_{d-2}a_d^{-1}a_0 & \ast  & \ast  &\cdots  & \ast\\
\vdots&\vdots&\vdots& &\vdots \\
-s_{1}a_d^{-1}a_0  & \ast  &\ast  & \cdots  & \ast\\
 \end{array}
\right)_{d\times d} \left(
                      \begin{array}{c}
                        1 \\
                        \beta_1 \\
                        \beta_2 \\
                        \vdots \\
                        \beta_{d-1} \\
                      \end{array}
                    \right)=\left(
                      \begin{array}{c}
                        \theta^{5i'+\ell-k} \\
                        0 \\
                        0 \\
                        \vdots \\
                        0 \\
                      \end{array}
                    \right)}\\
            &\Longleftrightarrow & \theta^{5i'+\ell-k}|S|=1~~\text{({ as  $\beta_i$ can be choosen freely from $\FF_e$ for $1\leq i\leq d-1$})}\\
\end{array}$$
By the definition of $\lg A \rg$, we may select an element $a:=\theta^{\delta_1}+s_1\alpha+s_2\alpha^2+\cdots+s_{d-1}\alpha^{d-1}\in\FF_{e^d}^*$ of order $e^d-1$, where this element corresponds to the element $\iota'\in \lg A\rg$ under our construction. Through computation of $a^2$, we find that $\sum_{i=1}^{d-1}s_{i}s_{d-i}$ appears as the coefficient of $\alpha^d$ in the expansion. This coefficient depends solely on the choice of the element $a$. Observe that $|\iota'|=\theta^{\delta_1}|S|+a_d^{-1}a_0\sum_{i=1}^{d-1}s_{i}s_{d-i}=1$. From this, we derive the equivalence:
$$t^{k-\ell}\iota' n'\in R \Longleftrightarrow
\theta^{5i'+\ell-k}=\theta^{\chi(a)}$$
where $\theta^{\chi(a)}:=\frac{\theta^{\delta_1}}{1-a_d^{-1}a_0(s_1s_{d-1}+s_2s_{d-2}+\cdots+s_{d-1}s_1)}$
denotes an element of $\FF_e^*$ whose value depends on the choice of $a$.
Let $$\Delta:=\{\chi(a^l)\mid \theta^{\chi(a^l)}=\frac{\theta^{\delta_l}}{1-a_{d}^{-1}a_{0}(s_{1,l}s_{d-1,l}+s_{2,l}s_{d-2,l}+\cdots+s_{d-1,l}s_{1,l})},  1\leq l\leq p\},$$ where $a^l=\theta^{\delta_l}+s_{1,l}\alpha+s_{2,l}\alpha^2+\cdots+s_{d-1,l}\alpha^{d-1}$.
Furthermore, the preceding arguments establish that:
 $$B_{k,1}\in \N(B_{\ell,0})\Longleftrightarrow\omega+k-\ell\equiv 0 \pmod 5~~~\text{for some $\omega\in\Delta$}.$$
 Let $\phi': B_{\mu, 0}\mapsto B_{\mu+1, 0}$, $B_{\vartheta, 1}\mapsto B_{\vartheta+1, 1}$ be a mapping from $V(X_{\mathcal{B'}})$ to $V(X_{\mathcal{B'}})$. It is not difficult to verify that the above mapping $\phi'$ is an automorphism of the quotient graph $X_{\mathcal{B'}}$. Consequently, we obtain the desired semiregular automorphism $\phi'$.

\vskip 3mm
{\it Step 1.2: Show the valency of $X_{\mathcal{B}'}$ is greater than $1$.}
\vskip 3mm
Applying the preceding analysis, we have that $B_{k,1}\in \N(B_{\ell,0})$ if and only if $Rt^k\lg A\rg\cap Rt^\ell N_1\neq \emptyset$. This condition is algebraically equivalent to $\theta^{5i+\ell-k}|S|=1$, for some ${\tiny \left(
  \begin{array}{cc}
    \theta^{\delta_1} & \gamma \\
    (\gamma')^T & S \\
  \end{array}
\right)\in\lg A\rg}$,
 where $\gamma=(s_1,s_2,\cdots,s_{d-1})$ and $\gamma'=(-s_{d-1}a_{d}^{-1}a_0,-s_{d-2}a_{d}^{-1}a_0,\cdots,\\-s_{1}a_{d}^{-1}a_0)$
 for some $s_i\in\FF_e$ and $1\leq i\leq d-1$. Derived from $Rt^\ell\in Rt^\ell\lg A \rg\cap Rt^\ell N_1$ that $B_{\ell,1}\in \N(B_{\ell,0})$ and $|S|\neq 0$.
To complete the proof of Step 1.2, it suffices to find an element  ${\tiny \left(
  \begin{array}{cc}
    \theta^{\delta_1} & \gamma \\
    (\gamma')^T & S \\
  \end{array}
\right)\in\lg A\rg}$ such that $|S|\notin (\FF_e^*)^5$. Rewrite $S=(\nu_1,\nu_2,\nu_3,\cdots,\nu_{d-1})$, where $\nu_i^T\in V(d-1,e)$. And let $S_{i'}:=(\nu_1,\nu_2,\cdots,\nu_{i'-1},\gamma'^T,\nu_{i'+1},\cdots,\nu_{d-1})$, where $2\leq i'\leq d-1$ and set $S_1:=(\gamma'^T,\nu_2,\nu_3,\cdots,\nu_{d-1})$.

We claim there exists an integer $m$ with $1\leq m\leq d-1$ such that $|S_m|\neq 0$. Suppose, for contradiction, that $|S_i|=0$ for all $1\leq i\leq d-1$. Recall that $|S|\neq 0$. Let $W_{i'}:=\lg \nu_1,\nu_2,\cdots,\nu_{i'-1},\nu_{i'+1},\cdots,\nu_{d-1} \rg$ and $W_1:=\lg \nu_2,\nu_3,\cdots,\nu_{d-2}\rg$ be subspaces of $V(d,e)$, where $2\leq i'\leq d-1$. Then $\gamma'^T\in (\bigcap_{i'=2}^{d-1}W_{i'}\bigcap W_1)$, i.e., $\gamma'=(0,0,\cdots,0)$, a contradiction to the definition of $\lg A \rg$, as $a_0\neq 0$.

Let $g:={\tiny \left(
  \begin{array}{cc}
    1 & \xi \\
    0 & E_{d-1} \\
  \end{array}
\right)}$, where $\xi=(0,0,0,\cdots,\underbrace{y}_{m},\cdots,0)\in V(d-1, e)$ and $y\in\FF_e^*$, that is, $\xi$ is a row vector whose $m$-th component is not equal to zero. And set ${\tiny \tau:=\left(
  \begin{array}{cc}
    \theta^{\delta_1} & \gamma \\
    (\gamma')^T & S \\
  \end{array}
\right)^g=\left(
  \begin{array}{cc}
    \theta^{\delta_1}-\xi \gamma'^T & (\theta^{\delta_1}-\xi \gamma'^T)\xi+(\gamma-\xi S) \\
    (\gamma')^T & \gamma'^T\xi + S \\
  \end{array}
\right).}$
Then $\tau$ is also a $(5,p)$-semiregular automorphism of $\Gamma$. Moreover, $|\gamma'^T\xi + S|=|S|+y|S_m|$. We can choose some $y$ such that $|\gamma'^T\xi + S|\notin(\FF_e^*)^5$.
Then, we can replace $\lg A\rg$ with $\lg \tau\rg$ to construct a new graph instead of $X_{\mathcal{B'}}$. Without loss of generality, this yields the desired result.
\vskip 3mm
{\it Step  1.3: Show a Hamilton cycle of  $X_{\B'}$ lifts.}
\vskip 3mm

Derived from Step 1.1 and Step 1.2 that any edge in  $X_{\mathcal{B'}}$ is contained in one of its Hamilton cycles, due to the action of the $(2,5)$-semiregular automorphism on $X_{\mathcal{B'}}$.
The arguments preceding Step 1.1 show that the valency of the graph $\Gamma$ is  $e^{d-1}$. Since $e$ is prime, $p=\frac{e^d-1}{e-1}$, $(d,e-1)=1$ and $5\di (e-1)$, it follows that $e\geqslant11$ prime and $d$ is odd, that is $d-1\geqslant2$ even and $e^{d-1}\geqslant 11^2$. However, the quotient graph $X_{\mathcal{B'}}$ has $10$ vertices, so there exists a block $B'\in N(B_{0,0})$ such that $d(B_{0,0},B')\geqslant2$. 
By Proposition~\ref{pro:4}, it follows that every  Hamilton cycle in $X_{\mathcal{B'}}$ containing the edge $\{B_{0,0},B'\}$ can be lifted to a Hamilton cycle in $X$, as required.

\vskip 3mm
{\bf Case 2}: $X[\mathcal{B}_1,\mathcal{B}_2]$ is the incidence graph of points and hyperplanes of $PG(d-1,e)$.
\vskip 3mm
Using the same notations and similar arguments as in Case 1, we conclude that the quotient graph $X_{\mathcal{B'}}$ induced by $\lg A\rg$, where $A$ is a $(5,p)$-semiregular automorphism,
has vertex set: $$V(X_{\mathcal{B'}})=\{B_{i,0}\mid 0\leq i\leq 4\}\cup \{B_{i,1}\mid 0\leq i\leq 4\},$$ where $$B_{i,0}=\{N_1Ht^is\mid s\in\lg A \rg\} \quad{\rm and}\quad B_{i,1}=\{N_2Ht^is\mid s\in\lg A \rg\}.$$
Now, we proceed with the proof in two steps.

\vskip 3mm
{\it Step 2.1: Show that $X_{\mathcal{B'}}\cong\K_{5,5}$.}
\vskip 3mm

Let  $w'=\{\lg (x_0, x_1, \cdots, x_{d-2}, 0)\rg_{\FF_{e}}\mid x_i\in\FF_{e}, 0\leq i\leq d-2\}$. Then $v\in w'$, where $v$ is defined as above.
We may set $\Gamma=\rm{B}(N, L, R', D')$, where $D'=R'L$ and $$R':=N_{w'}=\{ \left(\begin{array}{cc}
B' & \b'^T \\0 & \th^{5i'} \\
\end{array}\right)_{d\times d}\mid i'\in\ZZ, \beta'\in V(d-1, e), B'\in\GL(d-1, e)\}.$$
Let $I$ be a diagonal matrix defined as $I:=[\theta^{\epsilon},\theta^{-\epsilon-5e'},1,1,\cdots,1,\theta^{5e'}]
\in R'$, where $0\leq \epsilon\leq 4 $ is an integer and $e'$ is some integer. 
Then $(1,0,0,\cdots,0)I=(\theta^\epsilon,0,0,\cdots,0)$. And $I$ stabilizes $w'$ and maps the set $\{(\th^{5j},0,0,\cdots,0)\mid j\in\ZZ\}$ to the set $\{(\th^{5j+\epsilon},0,0,\cdots,0)\mid j\in\ZZ\}$ for any $0\leq \epsilon\leq 4$. Note that $\{(\th^{5j+\epsilon},0,0,\cdots,0)\mid j\in\ZZ\}$ is in the corresponding set $B_{\epsilon,0}$. Consequently, the quotient graph satisfies $X_{\mathcal{B'}}\cong \K_{5,5}$.
\vskip 3mm
{\it Step 2.2: Show that a Hamilton cycle in $X_{\mathcal{B'}}$ lifts to a Hamilton cycle in $X$.}
\vskip 3mm
By the result of Step 2.1, we have that the quotient graph $X_{\mathcal{B'}}\cong \K_{5,5}$. Then it is easy to deduce that any edge in the quotient graph $X_{\mathcal{B'}}$ is contained in one of its Hamilton cycles. On the other hand, since the valency of the graph $X[\mathcal{B}_1,\mathcal{B}_2]$ is $\frac{e^{d-1}-1}{e-1}$ and the quotient graph $X_{\mathcal{B'}}$ has $10$ vertices, there exists a block $B'\in \N(B_{0,0})$ such that $d(B_{0,0},B')\geqslant 2$. Consequently, by Proposition~\ref{pro:4}, any Hamilton cycle in $X_{\mathcal{B'}}$ containing the edge $\{B_{0,0},B'\}$ lifts to a Hamilton cycle in $X$, as desired.
\vskip 3mm
{\bf Case 3}: $X[\mathcal{B}_1,\mathcal{B}_2]\cong\K_{p,p}-p\K_{1,1}$ or $K_{p,p}$.
\vskip 3mm

For these two cases, it suffices to consider $X[\mathcal{B}_1,\mathcal{B}_2]\cong\K_{p,p}-p\K_{1,1}$, as $K_{p,p}-pK_{1,1}$ is a subgraph of $K_{p,p}$.
Therefore, we assume $X[\mathcal{B}_1,\mathcal{B}_2]\cong\K_{p,p}-p\K_{1,1}$. By the arguments preceding Case 1, this implies $\mathcal{B}_1=\mathcal{B}_2\in\{\mathcal{B}_v,\mathcal{B}_w\}$, and the graph $\Gamma=\rm{B}(N,L,L,D)$, where $D$ is a union of some double cosets of $L$ and $L$ in $N$, is a generating subgraph of $X$. In particular, $X[\mathcal{B}_1,\mathcal{B}_2]$ is a quotient graph of $\Gamma$. To complete the proof of Lemma \ref{quasi-pri-1}, it remains to show that $\Gamma$ admits a Hamilton cycle.

\vskip 3mm
{\it Claim 3.1: For any two adjacent vertices $\{B_a,B_b\}$ in the quotient graph $X[\mathcal{B}_1,\mathcal{B}_2]$ with $B_a\in\mathcal{B}_1$ and $B_b\in\mathcal{B}_2$, we have $X[B_a,B_b]\cong K_{5,5}$.}
\vskip 3mm

Note that $|N_{(A)}|=|N_{(B)}|$, for any $A,B\in \mathcal{B}_i$ (where $i=1,2$), as $N\trianglelefteq G$ and $N$ is $2$-transitive on either the projective points or hyperplanes. Moreover, $|N_{(A)}|=|N_{(C)}|$ holds for any $A\in \mathcal{B}_1$ and $C\in \mathcal{B}_2$, as $N\trianglelefteq G$ and there exists an automorphism $\tau\in G$ such that $\mathcal{B}_1^{\tau}=\mathcal{B}_2$. If there exists an element $g\in N_{(B_a)}\setminus N_{(B_b)}$, then $X[B_a,B_b]\cong K_{5,5}$ due to the action of $g$ on the set $B_b$. Now, it is sufficient to prove that $N_{(B_a)}\neq N_{(B_b)}$, where $B_a$ and $B_b$ are arbitrary two adjacent blocks in the quotient graph $X[\mathcal{B}_1,\mathcal{B}_2]$ such that $B_a\in\mathcal{B}_1$ and $B_b\in\mathcal{B}_2$.
On the contrary, we assume that $N_{(B_a)}= N_{(B_b)}$.
Since $N$ is transitive on the set $[N:L]$ and $N_{(B_a)}\trianglelefteq N$, we have that $N_{(B_a)}=N_{(B_b)}$
acts trivially on all vertices of $\Gamma$. It follows that $N_{B_a}\lesssim \S_5$ by the action of $N_{B_a}$ on $B_a$, a contradiction arisen as $\frac{|N|}{|N_{B_a}|}=p$ and $N\cong \SL(d,e)$. Hence the Claim 3.1 holds and so $X[B_a,B_b]\cong\K_{5,5}$.

However, since  $X[\mathcal{B}_1,\mathcal{B}_2]\cong K_{p,p}-pK_{1,1}$, one yields that $d(\Gamma)\geqslant 5(p-1)$, and hence $X$ contains a Hamilton cycle followed from Proposition \ref{pro:2}.

\vskip 3mm
{\bf Case 4}: $X[\mathcal{B}_1,\mathcal{B}_2]\cong pK_{1,1}$.
\vskip 3mm
No loss of generality, let $\mathcal{B}_1:=\{B_1,B_2,\cdots,B_p\}$ and $\mathcal{B}_2:=\{B_1',B_2',\cdots,B_p'\}$, where each $B_i$ is adjacent to $B_i'$ in the graph $X[\mathcal{B}_1,\mathcal{B}_2]$.
By the connectivity of $X$, there exists at least one edge within $\mathcal{B}_1$. Since $N$ is $2$-transitive on the set $\mathcal{B}_1$, we get that $B_i$ is adjacent to $B_j$ for any $1\leq i,j\leq p$.
By an argument analogous to Claim 3.1, we conclude that $X[B_i,B_j]\cong K_{5,5}$ for any $1\leq i,j\leq p$. Consequently, $d(X)\geq 5(p-1)$, and thus $X$ contains a Hamilton cycle by Proposition \ref{pro:2}.
\qed

\begin{lemma} \label{quasi-pri-2p} Under Hypothesis \ref{hypothesis1}, assume that $\soc(K)\cong\PSL(d,e)$ is faithful, imprimitive and quasiprimitive on an $N$-orbit of size $2p$, where $p=\frac{e^d-1}{e-1}$, both $d$ and $e$ are odd primes, $(d,e-1)=1$. Then $X$ contains a Hamilton cycle.
\end{lemma}
\demo
Suppose that $T:=\soc(K)\cong\PSL(d,e)\cong\SL(d,e)$, where $(d, e-1)=1$. Then $T\le N\le K\le \Aut(T)$
and  $V(X)$ is composed of $5$ $N$-blocks of size $2p$ each. Clearly, $\tau$ is transitive on these $5$ $N$-orbits, for any  element $\tau\in G\setminus K$ with $5$-power order.
To prove this lemma, it is sufficient to analyze the action of the group  $T\lg \tau\rg$ on the graph $X$. And, without loss of generality, we may assume that $N=T\cong \PSL(d,e)\cong \SL(d,e)$ and $G=N\lg\tau\rg$, where $|\tau|=5^s$ for some integer $s\geq 1$.

Let $\mathcal{B}=\{B_0,B_1,B_2,B_3,B_4\}$, where each $B_i$ is an $N$-orbit for $i\in\ZZ_5$. Since $\tau$ acts transitively on $\mathcal{B}$, the representations of $N$ on each $B_i\in\mathcal{B}$ are equivalent. The graph $X_{\mathcal{B}}$ contains a Hamilton cycle. Without loss of generality, we may assume $B_{i+1}\in \N(B_i)$ for all $i\in\ZZ_5$. 
Following the same method as in Lemma \ref{quasi-pri-1}, we decompose each  $B_{i}$ as $\bigcup_{j=1}^{p}B_{i,j}$ with $|B_{i,j}|=2$  for $i\in\ZZ_5$. Further, we may assume the action of $N$ on $\mathcal{B}_i:=\{B_{i,j}\di 1\leq j\leq p\}$ is equivalent to the action of $\PSL(d,e)$ on projective points (or hyperplanes). Since $|\mathcal{B}|=5$, the bipartite graph $X[\mathcal{B}_i,\mathcal{B}_{i+1}]$ cannot be isomorphic to the incidence graph (or non-incidence graph) of points and hyperplanes in $\PG(d-1,e)$.
We conclude that for each $i\in\ZZ_5$, the bipartite subgraph $X[\mathcal{B}_i,\mathcal{B}_{i+1}]$ is isomorphic to $(1)$ $K_{p,p}-pK_{1,1}$, $(2)$ $K_{p,p}$ and $(3)$ $pK_{1,1}$.

For the first two cases, applying the same method as in Claim $3.1$ of Lemma \ref{quasi-pri-1}, we conclude that $X[B_{i,a},B_{i+1,b}]\cong K_{2,2}$ holds for any adjacent pair $B_{i,a}\in \mathcal{B}_i$ and $B_{i+1,b}\in\mathcal{B}_{i+1}$ in $X[\mathcal{B}_i,\mathcal{B}_{i+1}]$, where $i\in\ZZ_5$.
Now, we get that $d(X)\geq 2(2p-2)\geqslant \frac {10p}3$, and so $X$ contains a Hamilton cycle, as Proposition \ref{pro:2}.

Next, consider the case where $X[\mathcal{B}_i,\mathcal{B}_{i+1}]\cong pK_{1,1}$ for all $i\in\ZZ_5$.
Without loss of generality, we may assume {\small $B_{i,j}\sim B_{i+1,j}$} for $i\in \ZZ_5$ and $1\leq j\leq p$, as {\small$B_{i,j}\sim B_{i+1,j}$} hlods if and only if $N_{B_{i,j}}=N_{B_{i+1,j}}$. By the connectivity of $X$, we have that either {\small $B_{0,1}\sim B_{k,t}$} or {\small $B_{0,1}\sim B_{0,t}$}  for some $k\in\ZZ_5\setminus\{0,1,4\}$ and $2\leq t\leq p$. Suppose {\small $B_{0,1}\sim B_{k,t}$} for some $k\in\ZZ_5\setminus\{0,1,4\}$ and $2\leq t\leq p$. Note that {\small$B_{i,j}\sim B_{i+1,j}$} if and only if $N_{B_{i,j}}=N_{B_{i+1,j}}$, for $i\in\ZZ_5$ and $1\leq j\leq p$. Then $N_{B_{0,1}}\neq N_{B_{k,t}}$ and so the graph $X[\mathcal{B}_0,\mathcal{B}_k]\cong K_{p,p}$ or $X[\mathcal{B}_0,\mathcal{B}_k]\cong K_{p,p}-pK_{1,1}$, as $N$ is $2$-transitive on each set $\mathcal{B}_i$, where $i\in\ZZ_5$. Now, we get a Hamilton cycle of $X$ by the first two cases and the action of $\tau$ on $\mathcal{B}$. Now suppose that {\small$B_{0,1}\sim B_{0,t}$} for some $2\leq t\leq p$. 
Since $N$ is $2$-transitive on the $p$-points, we have that {\small$B_{0,l}\sim B_{0,j}$} for any two distinct integers $1\leq l,j\leq p$. Following the proof technique of Lemma \ref{quasi-pri-1}, we can also get that $X[B_{0,l},B_{0,j}]\cong K_{2,2}$ for any $1\leq l,j\leq p$. Now, by the action of $\tau$ on the set $\mathcal{B}$, we get that $X[B_{i,l},B_{i,j}]\cong K_{2,2}$ for any $1\leq l,j\leq p$ and $i\in\ZZ_5$.

Now, we will get a Hamilton cycle of $X$. Note that in each vertex set $B_i$ induced subgraph, $X[B_{i,l},B_{i,j}]\cong K_{2,2}$ for any $i\in\ZZ_5$ and $1\leq l,j\leq p$.
For any set $B_i$ (where $i\in\ZZ_5$), the induced subgraph $X\lg B_i\rg$ contains a Hamilton path connecting any two distinct vertices of $B_i$.
Then we get five distinct $(2p-1)$-paths $P_i(\a_i,\b_i)$ in $X$ for $i\in\ZZ_5$, where each path $P_i(\a_i,\b_i)$ goes through all vertices in $B_i$ and its initial (resp. terminal) vertex  is $\a_i$ (resp. $\b_i$) such that $\{\a_i,\b_i\}=B_{i,1}$.
Let ${B}_{i,j}=\{x_{i,j},y_{i,j}\}$ for $i\in\ZZ_5$ and $1\leq j\leq p$. Given that $X[\mathcal{B}_k,\mathcal{B}_{k+1}]\cong pK_{1,1}$, we may assume (without loss of generality) the adjacency relations:   {\small$x_{k,j}\sim x_{k+1,j}$} and {\small$y_{k,j}\sim y_{k+1,j}$} for all $0\leq k\leq 3$ and $1\leq j\leq p$. 
We construct a Hamilton path of $X$ by concatenating the paths: 
{\small $P_0(x_{0,1}, y_{0,1} )\sim P_1(y_{1,1},x_{1,1})\sim P_2(x_{2,1},y_{2,1})\sim P_3(y_{3,1},x_{3,1})\sim P_4(x_{4,1},y_{4,1})$}, 
where each $P_i$ is the previous defined $(2p-1)$-path in $B_i$, and $\sim$ denotes adjacency between path endpoints. If {\small$y_{4,1}\sim x_{0,1}$} and {\small$x_{4,1}\sim y_{0,1}$}, 
then $X$ contains a Hamilton cycle: 
$${\small P_0(x_{0,1}, y_{0,1} )\sim P_1(y_{1,1},x_{1,1})\sim P_2(x_{2,1},y_{2,1})\sim P_3(y_{3,1},x_{3,1})\sim P_4(x_{4,1},y_{4,1})\sim P_0(x_{0,1},y_{0,1}).}$$ 
If {\small$y_{4,1}\sim y_{0,1}$} and {\small$x_{4,1}\sim x_{0,1}$,} then  
{\small$x_{i,j}\sim x_{i+1,j}$} and {\small$y_{i,j}\sim y_{i+1,j}$} for $i\in\ZZ_5$ and $1\leq j\leq p$, as the transitive action of $N$ on each vertex $B_i$. Then there exist $2p$ distinct $5$-cycles in $X$, which can be partitioned into two families: {\small$x_{0,j}\sim x_{1,j}\sim x_{2,j}\sim x_{3,j}\sim x_{4,j}\sim x_{0,j}$ and $y_{0,j}\sim y_{1,j}\sim y_{2,j}\sim y_{3,j}\sim y_{4,j}\sim y_{0,j}$} for each $1\leq j\leq p$. Then there exist $2p$ distinct $4$-paths, which decompose into two families:
$${\small P_{x,j}:=x_{0,j}\sim x_{1,j}\sim x_{2,j}\sim x_{3,j}\sim x_{4,j}}$$
and
$${\small P_{y,j}:=y_{0,j}\sim y_{1,j}\sim y_{2,j}\sim y_{3,j}\sim y_{4,j}}$$
 for each $1\leq j\leq p$.
 Note that for $t=0,4$, we have that {\small$x_{t,j}\sim x_{t,j'}$}, {\small$x_{t,j}\sim y_{t,j'}$} and {\small$y_{t,j}\sim y_{t,j'}$} for any $j\neq j'$ and $1\leq j,j'\leq p$. Define the reversed paths for each $1\leq j\leq p$:
$$ {\small P_{x,j}^{-}:=x_{4,j}\sim x_{3,j}\sim x_{2,j}\sim x_{1,j}\sim x_{0,j}}$$
and
$${\small P_{y,j}^{-}:=y_{4,j}\sim y_{3,j}\sim y_{2,j}\sim y_{1,j}\sim y_{0,j}.}$$
 Then, the following concatenation forms a Hamilton path in $X$ from $x_{0,1}$ to $y_{0,p}$:
$$ {\small
 P_{x,1}\sim P^{-}_{x,2}\sim P_{x,3}\sim\cdots\sim P_{x,p}\sim P_{y,1}^{-1}\sim P_{y,2}\sim P_{y,3}^{-}\cdots\sim P_{y,p}^{-1}.}$$ Since $x_{0,1}\sim y_{0,p}$, a Hamilton cycle of $X$ is desired at last.
\qed

\subsection{Two families of graphs where $K$ is unfaithful on a $r$-block}\label{K-unfaithful}
\begin{lemma}\label{unfaithful1} Under Hypothesis \ref{hypothesis1}, replacing $5$ by any prime $q$, assume $K$ acts unfaithfully on a $r$-block, where $r\in \{2, p, q\}$.
Then $X$ is hamiltonian.
\end{lemma}
\demo
Suppose that there exists a vertex $\b \in V(X)$ for which the group $K$ acts unfaithfully on its $N$-orbit $B_0:=\b^{N}$. Clearly, the normal quotient graph $X_{\B}$ is a $G/K$ vertex-transitive connected graph of order $m$, where $m\in\{2q, pq,2p\}$.
By Proposition \ref{pq-H-cycle}, $X_\B$ admits a Hamilton cycle unless it is the Petersen graph. For $X$, applying Proposition~\ref{pro:3}, one yields that there exists a $(m, r)$-semiregular automorphism $\rho\in G$ whose orbits coincide with the $N$-orbits, where $r:=\frac{2pq}{m}$.
Then $G_{B_0}=HK$, where $H:=G_{\b}$, $N\leq K$ and $\rho\in HK$. We now proceed to get a Hamilton cycle in $X$.

First, assume that $X_{\B}$ is not isomorphic to the Petersen graph.
Then $X_\B$ contains a Hamilton cycle $C_{\mathcal{B}}$ of the form:
$B_0\sim B_1\sim B_2\sim \cdots\sim  B_{m-1}\sim B_0,$
where $B_i\in\mathcal{B}$ for $i\in\ZZ_m$.
Note that $K_{(B_0)}\neq1$.
Considering the action of $K_{(B_0)}$ on each set $B_i$, as $K_{(B_0)}\unlhd K$ and $|B_i|=r$ is a prime for $i\in\ZZ_m$, we can conclude that there must exist two adjacent $N$-orbits, say $B_j$ and $B_{j+1}$ such that $d(B_j, B_{j+1})=r>1$, for some $j\in\ZZ_m$.
By Proposition~\ref{pro:4}, $C_{\mathcal{B}}$ can be lifted to a Hamilton cycle of $X$.
So in what follows, we assume that $X_{\B}$ is isomorphic to the Petersen graph.

\begin{figure}\centering
{\small \begin{tikzpicture}[every node/.style={draw, circle}]
    % First pentagon

    \node (B0) at (-2,4.5) {$B_0$};
    \node (B1) at (1,2.5) {$B_1$};
    \node (B2) at (0,-1) {$B_2$};
    \node (B3) at (-4,-1) {$B_3$};
    \node (B4) at (-5,2.5) {$B_4$};
    \node (B5) at (-2,3) {$B_5$};
    \node (B6) at (-0.5,2) {$B_6$};
    \node (B7) at (-1,0) {$B_7$};
    \node (B8) at (-3,0) {$B_8$};
    \node (B9) at (-3.5,2) {$B_9$};

    \draw (B0) -- (B1) -- (B2) -- (B3) -- (B4) -- (B0);
    \draw (B0) -- (B5);
    \draw (B1) -- (B6);
    \draw (B2) -- (B7);
    \draw (B3) -- (B8);
    \draw (B4) -- (B9);
    \draw (B5) -- (B7) -- (B9) -- (B6) -- (B8) -- (B5);

    % Second pentagon
    \node (B'0) at (7,4.5) {$B_0$};
    \node (B'1) at (10,2.5) {$B_4$};
    \node (B'2) at (9,-1) {$B_2$};
    \node (B'3) at (5,-1) {$B_1$};
    \node (B'4) at (4,2.5) {$B_3$};
    \node (B'5) at (7,3) {$B_9$};
    \node (B'6) at (8.5,2) {$B_7$};
    \node (B'7) at (8,0) {$B_5$};
    \node (B'8) at (6,0) {$B_8$};
    \node (B'9) at (5.5,2) {$B_6$};

    \draw (B'0) -- (B'1) -- (B'2) -- (B'3) -- (B'4) -- (B'0);
    \draw (B'0) -- (B'5);
    \draw (B'1) -- (B'6);
    \draw (B'2) -- (B'7);
    \draw (B'3) -- (B'8);
    \draw (B'4) -- (B'9);
    \draw (B'5) -- (B'7) -- (B'9) -- (B'6) -- (B'8) -- (B'5);
    \node[align=center, above=1cm of B0, anchor=north, draw=none]  {\small{$\Gamma_1$}};
\node[align=center, above=1cm of B'0, anchor=north, draw=none]  {\small{$\Gamma_2$}};
\end{tikzpicture}}
  \caption{ $\Gamma_1$: Petersen graph\,\,\,\,\,\,\,\,\,\,\,\,\,\,\,\,\,\,\,\,\,\,\,\,\,\,\,\,\,\,\,\, $\Gamma_2$: a subgraph of the complement graph of $\Gamma_1$}\label{coPetersengraph}
  \end{figure}
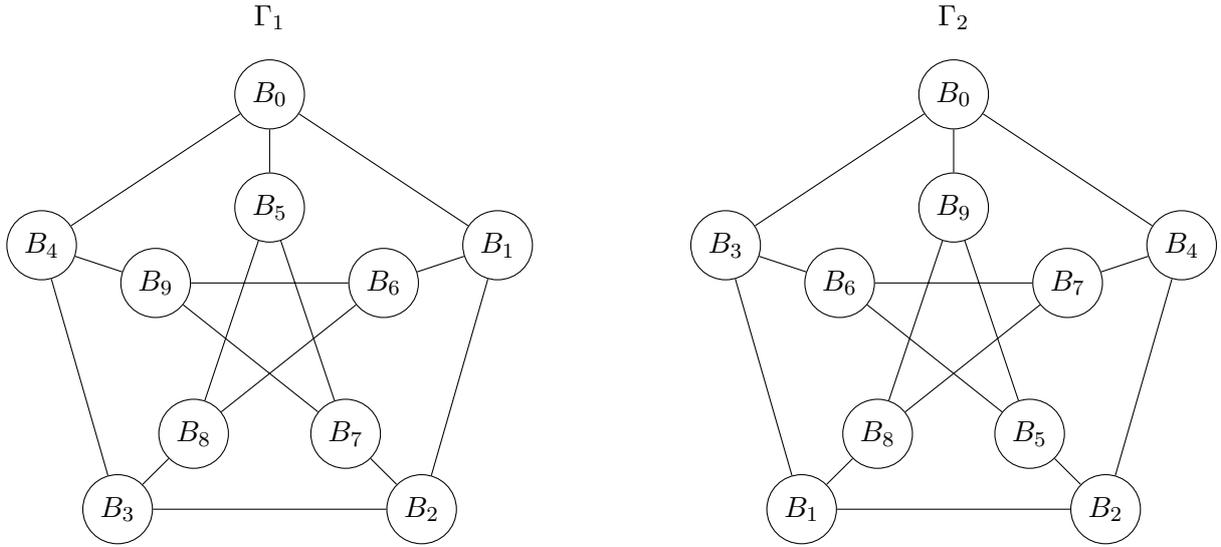

 For the case where $X_{\B}$ is isomorphic to the Petersen graph, we have $m=10$ and $r\geqslant7$. Without loss of generality, set $\mathcal{B}=\{B_0,B_1,\cdots,B_{9}\}$, where $B_i:=\{b_{i,0},b_{i,1},\cdots,b_{i,r-1}\}$, for $i\in\ZZ_{10}$. Then $X_{\B}$ is presented in Figure \ref{coPetersengraph}.
By analyzing  Figure \ref{coPetersengraph} carefully, for any two non-adjacent vertices $B_x,B_y\in V(X_{\mathcal{B}})$, there exists a Hamilton path $W_{xy}$ in $X_{\mathcal{B}}$ connecting $B_x$ to $B_y$, where $x,y\in\ZZ_{10}$. Representative examples include:
$${\small\begin{array}{ll}
&W_{02}:\,B_0\rightarrow B_4\rightarrow B_3\rightarrow B_8\rightarrow B_5\rightarrow B_7\rightarrow B_9\rightarrow B_6\rightarrow B_1\rightarrow B_2,\\
&W_{06}:\,B_0\rightarrow B_1\rightarrow B_2\rightarrow B_3\rightarrow B_4\rightarrow B_9\rightarrow B_7\rightarrow B_5\rightarrow B_8\rightarrow B_6.\\
\end{array}}$$
Given $K_{(B_0)}\neq 1$, now we consider the action of $K_{(B_0)}$ on each set $B_i$ for $i\in\ZZ_{10}$. Since $|B_i|=r$ is a prime, and $K_{(B_0)}\neq 1$, there exist two adjacent vertices $B_l,B_{l'}$ in $W_{02}$ such that $X[B_l,B_{l'}]\cong K_{r,r}$.
It is a fundamental property that any vertex-transitive subgroup of the Petersen graph's automorphism group has at most two orbits when acting on the edge set. Thus, the group $G/K$ acts on the edge set of $X_{\mathcal{B}}$ with at most two orbits, which are:
$$ {\small
\begin{array}{lll}
\mathcal{E}&:=&\{\{B_0,B_4\}, \{B_4,B_3\}, \{B_3,B_2\}, \{B_2,B_1\}, \{B_1,B_0\}, \{B_5,B_8\}, \\
&&\{B_8,B_6\}, \{B_6,B_9\}, \{B_9,B_7\}, \{B_7,B_5\}\}
\end{array}}$$
 and
$$  {\small
\begin{array}{lll}
\mathcal{E}'&:=&\{\{B_0,B_5\}, \{B_4,B_9\}, \{B_3,B_8\}, \{B_2,B_7\}, \{B_1,B_6\}\}.
\end{array}}$$
These arguments establish the fact that
$X[B, C]\cong\K_{r, r}$ for either $\{B,C\}\in\mathcal{E}$ or $\{B,C\}\in\mathcal{E}'.$

If $X[B, C]\cong\K_{r, r}$ and $X[B',C']\cong rK_{1,1}$, for $\{B,C\}\in\mathcal{E}$ and $\{B',C'\}\in \mathcal{E}'$, then without loss of generality, define the adjacency relation between two blocks $\{B',C'\}\in \mathcal{E}'$ in the following:

$${\small\left\{
   \begin{array}{ll}
     b_{6,0}\sim b_{1,1}, b_{6,2}\sim b_{1,r-1},&  \\
     b_{6,\jmath+1}\sim b_{1,\jmath}, & \hbox{for $\jmath=0$ or $2\leq \jmath\leq p-2$,} \\
     b_{\imath,\jmath}\sim b_{5+\imath,\jmath+1}, & \hbox{for $\imath\in\ZZ_5$ and other cases.} \\
   \end{array}
 \right.}$$
Then we can get two different types of walks $W_{\a}$, $W_{\b,\jmath'}$, where
$${\small\begin{array}{ll}
W_{\a}:=& b_{0,0}\sim b_{4,0}\sim b_{3,0} \sim b_{8,1}\sim b_{5,0}\sim b_{7,0} \sim b_{9,1}\sim b_{6,0}\sim b_{1,1}\sim b_{2,0}\sim b_{7,1}\sim b_{5,1}\sim \\
&b_{8,0}\sim b_{6,1}\sim b_{9,0}\sim b_{4,1}\sim b_{3,1}\sim b_{2,2}\sim b_{1,0}\sim b_{0,1}\sim b_{5,2}\sim b_{8,2}\sim b_{6,2}\sim b_{9,2}\sim \\
& b_{7,2}\sim b_{2,1}\sim b_{3,2}\sim b_{4,2}\sim b_{0,2}\sim b_{1,2},
\end{array}}$$
and
$${\small\begin{array}{ll}
W_{\b,\jmath'}:=& b_{6,\jmath'}\sim b_{8,\jmath'}\sim b_{5,\jmath'} \sim b_{7,\jmath'}\sim b_{9,\jmath'+1}\sim b_{4,\jmath'} \sim b_{3,\jmath'}\sim b_{2,\jmath'+1}\sim b_{1,\jmath'}\sim b_{0,\jmath'}\sim b_{5,\jmath'+1}\\
&\sim b_{8,\jmath'+1}\sim b_{6,\jmath'+1}\sim b_{9,\jmath'}\sim b_{7,\jmath'+1}\sim b_{2,\jmath'}\sim b_{3,\jmath'+1}\sim b_{4,\jmath'+1}\sim b_{0,\jmath'+1}\sim b_{1,\jmath'+1},
\end{array}}$$
for $3\leq\j'\leq p-2$.
Then the Hamilton cycle
{\small$W_{\a}\sim W_{\b,3}\sim W_{\b,5}\sim\cdots\sim W_{\b,p-2}\sim b_{0,0}$} of $X$ is
obtained.

If $X[B, C]\cong r K_{1,1}$ and $X[B',C']\cong K_{r,r}$, for $\{B,C\}\in\mathcal{E}$ and $\{B',C'\}\in \mathcal{E}'$, then without loss of generality, define the adjacency relation between two blocks $\{B',C'\}\in \mathcal{E}'$ in the following:
$${\small\left\{
   \begin{array}{ll}
     b_{\imath,\jmath}\sim b_{\imath+1,\jmath+1}, & \hbox{for $\imath\in\ZZ_5$ and $\jmath\in\ZZ_r$,} \\
     b_{5+\imath,\jmath}\sim b_{5+(\imath+2),\jmath+1}, & \hbox{for $\imath\in\ZZ_5$ and $\jmath\in\ZZ_r$.} \\
   \end{array}
 \right.}$$
Then we can get three different types of walks $W_{\a}$, $W_{\b}$ and $W_{\gamma,\jmath'}$, where
$${\small\begin{array}{ll}
W_{\a}:=& b_{0,0}\sim b_{5,2}\sim b_{8,1} \sim b_{3,0}\sim b_{4,1}\sim b_{9,3} \sim b_{7,2}\sim b_{2,2}\sim b_{1,1}\sim b_{6,0}\sim b_{1,2}\sim b_{6,2}\sim \\
&b_{9,1}\sim b_{4,0}\sim b_{0,1}\sim b_{5,1}\sim b_{8,0}\sim b_{3,2}\sim b_{2,1}\sim b_{7,0}\sim b_{2,0}\sim b_{7,1}\sim b_{9,2}\sim b_{4,2}\sim \\
& b_{3,1}\sim b_{8,2}\sim b_{6,1}\sim b_{1,4}\sim b_{0,3}\sim b_{5,0},
\end{array}}$$
$${\small\begin{array}{ll}
W_{\b}:=& b_{0,2}\sim b_{5,4}\sim b_{8,3} \sim b_{3,3}\sim b_{4,4}\sim b_{9,4} \sim b_{7,3}\sim b_{2,4}\sim b_{1,3}\sim b_{6,4}\sim b_{1,0}\sim b_{6,3}\sim \\
&b_{8,4}\sim b_{3,4}\sim b_{2,3}\sim b_{7,4}\sim b_{9,5}\sim b_{4,3}\sim b_{0,4}\sim b_{5,3}
\end{array}}$$
and
$${\small\begin{array}{ll}
W_{\gamma,\jmath'}:=& b_{0,\jmath'}\sim b_{5,\jmath'+1}\sim b_{8,\jmath'} \sim b_{3,\jmath'}\sim b_{4,\jmath'+1}\sim b_{9,\jmath'+1} \sim b_{7,\jmath'}\sim b_{2,\jmath'+1}\sim b_{1,\jmath'}\sim b_{6,\jmath'+1}\sim \\
&b_{1,\jmath'+1}\sim b_{6,\jmath'}\sim b_{8,\jmath'+1}\sim b_{3,\jmath'+1}\sim b_{2,\jmath'}\sim b_{7,\jmath'+1}\sim b_{9,\jmath'+2}\sim b_{4,\jmath'}\sim b_{0,\jmath'+1}\sim b_{5,\jmath'},
\end{array}}$$
for $5\leq\j'\leq p-2$.
Then the Hamilton cycle
{\small$W_{\a}\sim W_{\b}\sim W_{\gamma,5}\sim W_{\gamma,7}\sim W_{\gamma,9} \cdots\sim W_{\gamma,p-2}\sim b_{0,0}$} of $X$ is
obtained.
\qed

\begin{lemma}\label{unfaithful2}
 Under Hypothesis \ref{hypothesis1}, replacing $5$ by any prime $q$, assume $N$ induces $r$-blocks, where $r\in\{2p,2q,pq\}$.
 Let $N_0$ be a minimal normal subgroup of $G$ such that $N_0$ induce the same blocks with $N$. If $N_0$ acts unfaithfully on an $N$-block, then $X[A,B]\cong K_{r,r}$ for any two adjacent $N$-blocks $A$ and $B$ in the graph $X_{\mathcal{B}}$. Therefore, $X$ is hamiltonian.
\end{lemma}
\demo
Let $s=\frac{2pq}{r}$, where $r\in\{2p,2q,pq\}$, and $\B=\{\a^N\di \a\in V(X)\}$ the set of all $N$-orbits of size $r$ each.
Since $G$ is a minimal transitive group on $V(X)$, it follows that $G/K$ is transitive on $V(X_\B)$.
Then there exists a subgroup $\lg\t\rg K/K$ of $G/K$ such that $\lg\t\rg K/K$ is a $s$-group and acts transitively on $V(X_\B)$. Consequently, $K\lg\t\rg$ acts transitively on $V(X)$.
By the minimality of $G$, we get $G=K\lg\t\rg$, where $\o(\tau)=s^{e'}$ for some integer $e'$ and $\lg\tau\rg$ acts transitively on $\B$.

Suppose that $N_0$ acts unfaithfully on an $N$-orbit $A=\a^N$. 
Let ${N_0}_{(A)}$ be the kernel of $N_0$ acting on the set $A$.
Noting that ${N_0}^A\cong N_0/{N_0}_{(A)}$ and $r\in\{2p,2q,pq\}$, it follows from Proposition~\ref{squarefree} that ${N_0}^A$ is a non-abelian simple group.

Let $B$ be a vertex adjacent to $A$ in the graph $X_{\mathcal{B}}$.
Then $({N_0}_{(B)}{N_0}_{(A)})/{N_0}_{(A)}$ is a normal subgroup of $N_0/{N_0}_{(A)}$. Then either ${N_0}_{(A)}={N_0}_{(B)}$ or $N_0={N_0}_{(A)}{N_0}_{(B)}$, as ${N_0}^A$ is a non-abelian simple group.
Suppose for a contradiction that ${N_0}_{(A)}={N_0}_{(B)}$.
Then ${N_0}_{(A)}={N_0}_{(C)}$, for any $C\in\mathcal{B}$, as the connectivity of the graph $X_{\mathcal{B}}$ and $\tau$ is transitive on $\mathcal{B}$. And then ${N_0}_{(A)}=1$, a contradiction.
Therefore, $N_0={N_0}_{(A)}{N_0}_{(B)}$.
Next, consider the transitive action of $N_0$ on $B$. Then we get that $X[A,B]\cong K_{r,r}$, as desired.
\qed

\section{Cases $r\in\{2,5,p\}$}\label{p-blocks}
In this section,  we address the case where   $K$ induces blocks of length $r$ under Hypothesis \ref{hypothesis1}, with  $r\in \{ 2, 5, p\}$.  The subcases $r=p$ and $r\in \{2, 5\}$ are discussed in Subsections 4.1 and 4.2, respectively.
\subsection{$r=p$}
First, we consider the case where $m=|V(X_\B)|=10$, i.e., $r=p$. Fix  $\a\in B\in\B$.

\begin{theorem}\label{p-blocks}
The graph $X$ is  Hamiltonian if $r=p$.
\end{theorem}
\pf Recall that $K$ is the kernel of $G$ on $\mathcal{B}$. If $K$
 acts unfaithfully on any of its orbits, then by Lemma \ref{unfaithful1}, $X$ contains a Hamilton cycle, as desired. Thus, in what follows, we assume that $K$ acts faithfully on each $K$-block, i.e., $\ZZ_p\lesssim K\leq \S_p$ and $K$ is a primitive permutation group of degree $p$. Let $\ZZ_p\cong P=\lg a\rg\leq K$.
Then $P\in\Syl_p(K)$, and by the Proposition \ref{the:main5}, we conclude that $G=\N_G(P)K$ and $\N_G(P)$ acts transitively on $V(X)$.
By the minimality of $G$, it follows that $G=\N_G(P)$, i.e., $P\lhd G$ and $P\lhd K$. Employing the classification of primitive permutation groups (see, e.g., \cite{LPS88}), we deduce that $K$ is of affine type and $K\leq\AGL(1, p)$. Thus, $K\cong\ZZ_p:\ZZ_l$ for some $l\di (p-1)$.
By the minimality of $G$ again, $G/K$ is a minimal transitive group of degree 10.
Lemma~\ref{deg10} implies that for some positive integer $\kappa$,
$$G/K\in\{\A_5,\ZZ_5:\ZZ_4,\ZZ_{10},\D_{10},\ZZ_5^2:\ZZ_{2^\kappa}, \ZZ_2^4:\ZZ_5\}.$$
Then to prove Theorem \ref{p-blocks}, we apply  Lemmas~\ref{H-Z5Z4},~\ref{A5},~\ref{Z10} and~\ref{Z52} (to be established later).
\qed

\begin{lemma}\label{cga}\label{Z5Z4}
Using the notations of Theorem \ref{p-blocks} and its proof, if $G/K\cong\ZZ_5:\ZZ_4$, then there exist elements $a,b,c,\sigma\in G$ such that$:$
\begin{itemize}
    \item[{\rm(1)}] $G=\lg a\rg\times \lg b \rg:\lg \sigma \rg$ and $K=\lg a \rg:\lg c \rg$, where $a^p=b^5=1$, $b^\sigma=b^2$, $c=\sigma^{4}$, $|\sigma|=2^{k+2}$ and $|c|=l=2^k$ for some integer $k\geq 0;$
    \item[{\rm(2)}]and either:
    \begin{enumerate}
     \item[{\rm(a)}] $\C_G(P)\cong\ZZ_{5p},$ or
      \item[{\rm(b)}] $l=1$ and $\C_G(P)\cong\ZZ_{5p}:\ZZ_2=\ZZ_p\times(\ZZ_5:\ZZ_2)$ or $\ZZ_{5p}:\ZZ_4=\ZZ_p\times(\ZZ_5:\ZZ_4).$
      \end{enumerate}
  \end{itemize}
\end{lemma}

\demo
(1)
First, we claim that $5p\di|\C_G(P)|$.
Note that $\lg a\rg=P\lhd G$, so $G/\C_G(P)\lesssim\Aut(P)$ is cyclic. Now, observe that $\C_G(P)K/K\lhd G/K$, and then $(G/K)/(\C_G(P)K/K)$ is isomorphic to $ G/\C_G(P)K\cong(G/\C_G(P))/(\C_G(P)K/\C_G(P))$. Thus, $(G/K)/(\C_G(P)K/K)$ is cyclic, as $G/\C_G(P)$ is cyclic. Note that $G/K\cong\ZZ_5:\ZZ_4$ has a unique minimal normal subgroup, isomorphic to $\ZZ_5$. Consequently, $(G/K)/(\C_G(P)K/K)$ must be isomorphic to $\ZZ_4$, $\ZZ_2$ or $1$. It implies that $\C_G(P)K/K$ is isomorphic to $\ZZ_5$, $\ZZ_5:\ZZ_2$ or $\ZZ_5:\ZZ_4$, respectively. Then there exists a subgroup $H\leq\C_G(P)K$ such that $H/K\cong\ZZ_5$ and $H\cong K.\ZZ_5\cong(\ZZ_p:\ZZ_l).\ZZ_5$ for $l\di (p-1)$.
On the other hand, since $\ZZ_p:\ZZ_l\cong K\lesssim\AGL(1, p)$, we conclude that  $\C_G(P)\cap K=P$ and $|H|=5pl\mid |\C_G(P)K|$. Then $|\C_G(P)K|=\frac{|\C_G(P)||K|}{|\C_G(P)\cap K|}=
|\C_G(P)|l$. Thus, $5p\di|\C_G(P)|$. This proves the claim.

Now, take an element $b\in\C_G(P)$ of order $5$. Then $\lg ab \rg=\lg a \rg\times\lg b \rg\leq \C_G(P)$ and $K\lg b\rg=K:\lg b \rg\leq \C_G(P)K$. Since $(K\lg b\rg)/K\leq G/K$ and $(K\lg b\rg)/K\cong \ZZ_5$, while $G/K\cong \ZZ_5:\ZZ_4$, we conclude that $K\lg b \rg\trianglelefteq G$. Since $\C_G(P)\cap K\lg b\rg=\lg ab \rg$ and $\C_G(P)\trianglelefteq G$, it follows that $\lg ab \rg\trianglelefteq G$. Note that $P\leq \lg ab \rg$. 
Therefore, each $K$-orbit is contained in a $\lg ab \rg$-orbit on $V(X)$. Since $b\notin K$ and $G/K\cong \ZZ_5:\ZZ_4$ is transitive on the complete $K$-block system. We get that $\lg ab \rg$ induces totally two blocks in $V(X)$, say $\Delta_1$ and $\Delta_2$, where $V(X)=\Delta_1\cup\Delta_2$ and $|\Delta_1|=|\Delta_2|=5p$. By the transitivity of $G$ on $V(X)$, there exists an element $\sigma\in G$ of $2$-power order such that $\Delta_1^\sigma=\Delta_2$. By the minimal transitivity of $G$, we may assume $G=\lg ab \rg:\lg \sigma \rg\cong \ZZ_{5p}:\ZZ_{2^{k+2}}$ for some integer $k\geq 0$, as $G/K\cong \ZZ_5:\ZZ_4$. Moreover, we can write $K=\lg a \rg:\lg c\rg\trianglelefteq G$, where $c\in G$ is an element of order $l$.

 We now claim that $b^\sigma=b^2$ and $c\in\lg \sigma \rg$. Clearly, $\lg b \rg\trianglelefteq G$, as $\lg b \rg$ is a character subgroup of $\lg ab \rg$. Without loss of generality, we may assume that $b^{\sigma}\in\{ b,b^2\}$. Suppose, for contradiction, that $b^\sigma=b$. Then $\lg b \rg:\lg \sigma \rg=\lg b\rg\times \lg \sigma \rg$, which is an abelian group. However, this contradicts the fact that $G/K\cong \ZZ_5:\ZZ_4$ is not abelian. So $b^\sigma=b^2$. If $5\mid |K|$, then $\lg b\rg\leq K$, as $|G|=5p2^{k+2}$ and $\lg b \rg\trianglelefteq G$ is the only subgroup of $G$ with order $5$. Since $\C_G(P)\cap K=P$, we get that $\lg b\rg\nleq K$. Consequently, $|K|\mid p2^{k+2}$. Therefore, $l\mid 2^{k+2}$, that is, $l=2^s$ for some integer $s\geq 0$. Since all Sylow $2$-subgroups of $G$ are conjugate in $G$ and $K\trianglelefteq G$, we may assume without loss of generality that $c=\sigma^{2^{k+2-s}}\in \lg \sigma \rg$, where $l=2^s$. 
 Since $\lg \sigma\rg/(K\cap \lg \sigma \rg)\cong(\lg \sigma\rg K)/K\cong \ZZ_4$ and $K\cap \lg \sigma \rg=\lg c \rg$,  we get that $|\lg c \rg|=\frac{2^{k+2}}{4}=2^k$. Given that $|\sigma|=2^{k+2}$ and $c\in\lg \sigma \rg$, we may assume without loss of generality that $c=\sigma^4$. Then we get that $s=k$, as $|\lg c\rg|=2^s=2^k$, the results are desired.
\vskip 3mm

(2) Observe that $\C_G(P)/P\cong\C_G(P)/(\C_G(P)\cap K)\cong \C_G(P)K/K$ which is isomorphic to $\ZZ_5$, $\ZZ_5:\ZZ_2$ or $\ZZ_5:\ZZ_4$. Note that $\lg c\rg=\lg\s^4\rg\cap\C_G(P)=1$. If $\ZZ_{5p}:\ZZ_2\lesssim \C_G(P)$, then $\lg \s^{2^{k+1}}\rg\leq \C_G(P)$. If $k\geq 1$, then $(\s^{4})^{2^{k-1}}\in \C_G(P)$, which contradicts the fact that $\lg c\rg=\lg\s^{4}\rg\cap\C_G(P)=1$.
We conclude that either $\C_G(P)\cong\ZZ_{5p}$; or $k=0$ and one of the following holds: $\C_G(P)\cong\ZZ_{5p}:\ZZ_2=\ZZ_p\times(\ZZ_5:\ZZ_2)$ or $\ZZ_{5p}:\ZZ_4=\ZZ_p\times(\ZZ_5:\ZZ_4)$, as desired. \qed

\begin{lemma}\label{H-Z5Z4}
Using the notations of Theorem \ref{p-blocks}, if $G/K\cong\ZZ_5:\ZZ_4$, then
$X$ contains a  Hamilton cycle.
\end{lemma}
\vskip 3mm
\demo Following the notations of Lemma \ref{Z5Z4}, let $\Omega:=\{\Delta_0,\Delta_1\}$ denote the set of $\lg ab \rg$-orbits on $V(X)$, where $|\Delta_i|=5p$ for $i=0,1$. We observe that each $K$-orbit is contained in a $\lg ab\rg$-orbit, and the set $\B=\{B_0, B_1, \cdots, B_9\}$, consisting of all $K$-orbits on $V(X)$, forms a complete block system for the action of $G$ on $V(X)$. Without loss of generality, we may assume that $B_j\subset \Delta_0$ and $B_{5+j}\subset \Delta_1$ for each $j\in\ZZ_5$.

Let $M:=\lg\s^2\rg$, which is a point stabilizer of $\lg\s\rg$ acting on $\Omega$. Since $M$ is a $2$-group and $|\Delta_0|=|\Delta_1|=5p$, it follows that $M$ is also a point stabilizer of $G$ acting on $V(X)$. We now represent the graph $X$ as the coset graph of group $G$ with respect to the subgroup $M$. Under this construction, we get that $V(X)=[G:M]$ and $G$ acts on $V(X)$ by right multiplication. The undirected nature of $X$ implies the following adjacency property: if
$M\sim Mg$ for $g\in G$, then $M\sim Mg^{-1}$.
 Note that $ab$ is a semiregular automorphism of order $5p$ with exactly two orbits $\Delta_0$ and $\Delta_1$ on $V(X)$. Let $Y_0=X\lg\Delta_0\rg$ and $Y_1=X\lg\Delta_1\rg$ be the induced subgraphs. Then $Y_0^\s=Y_1$ and $V(Y_i)=\{M\s^ia^xb^y\di x\in\ZZ_p, y\in\ZZ_5\}$ for $i=0, 1$.

Using the coset perspective, we relabel $B_{5i+j}$ as $M\s^i\lg a\rg b^j$ for $i=0,1$ and $j\in\ZZ_5$. Consider the quotient graph $X_\B$. Then $G/K$ is a vertex-transitive automorphism subgroup of $X_\B$ with a point stabilizer $MK/K\cong M/(K\cap M)=M/\lg c\rg=\lg \s^2\rg/\lg \s^4\rg\cong\ZZ_2$. Note that the edge set of $X_{\mathcal{B}}$ is the union of some orbitals of $G/K$ acting on $\mathcal{B}$. By analyzing the orbits of $MK/K$ on $\B=\{M\s^i\lg a\rg b^j|i=0, 1, j=0, 1, 2, 3, 4\}$, we conclude that except for the case $X_\B\cong K_{10}$, 
the normal quotient graph $X_\B$ is isomorphic to the Petersen graph, its complement, or the graph $\K_{5,5}- 5\K_{1,1}$ (also can be checked by Magma).
We observe that all three graphs (the Petersen graph, its complement, and $K_{5,5}-5K_{1,1}$) are all subgraphs of $K_{10}$. Therefore, it suffices to restrict our analysis to these three cases in subsequent arguments.
\vskip 3mm
 \textbf{Claim 1:} If $X_\B$ is the Petersen graph, then $X$ contains a Hamilton cycle.
\vskip 3mm
Suppose $X_\B$ is isomorphic to the Petersen graph.
Then $\N(B_0)=\{ B_1,B_4,B_5\}$, up to isomorphic, as illustrated in Figure\ref{coPetersengraph}. Consequently, there exist $i_1, i_2\in \ZZ_p$ such that: $$\{Ma^{i_1}b, Ma^{-i_1}b^{-1}, M\s a^{i_2}\}\subseteq\N(M).$$
Moreover, Lemma \ref{Z5Z4} establish that  $\C_G(P)\cong\ZZ_{5p}, \ZZ_{5p}:\ZZ_2=\ZZ_p\times(\ZZ_5:\ZZ_2)$ or $\ZZ_{5p}:\ZZ_4=\ZZ_p\times(\ZZ_5:\ZZ_4)$.
\vskip 3mm
 {\it Case 1:  $\C_G(P)\cong\ZZ_{5p}$}
\vskip 3mm

  Followed from the fact that $\ZZ_{5p}\cong\lg a\rg\times\lg b\rg\leq\C_G(P)\cong\ZZ_{5p}$, it follows that $\lg \sigma \rg\cap \C_G(P)=1$. Now, let $a^{\sigma}=a^t$ for some integer $t$.
Then $a^{t^{2^{k+1}}}=a^{\sigma^{2^{k+1}}}$ and $t^{2^{k+1}}\not\equiv1 (\bmod\, p)$, as $a\neq a^{\sigma^{2^{k+1}}}$. However, since $o(\s)=2^{k+2}$, we derive that $a=a^{\s^{2^{k+2}}}=a^{t^{2^{k+2}}}$ and $(t^{2^{k+1}})^2=t^{2^{k+2}}\equiv1 (\bmod\, p)$. Consequently, $t^{2^{k+1}}\equiv-1 (\bmod\, p)$ and $(p, t)=1$.
Since $\s\not\in\C_G(P)$, the `NC'-Theorem applied to $P:\lg\s\rg$ implies that $\s\in\Aut(P)\cong\ZZ_{p-1}$ and $2^{k+2}\mid p-1$. Consequently, the group $G$ admits the following presentation:
$$G=\lg a,b,\s\di a^p=b^5=\s^{2^{k+2}}=1,[a,b]=1,a^\s=a^t,b^\s=b^2, t^{2^{k+1}}\equiv-1 (\bmod\, p)\rg. $$
By the vertex-transitivity of $G$ on $X$, we derive that $\{Ma^{i_1}b\s, Ma^{-i_1}b^{-1}\s, M\s a^{i_2}\s\}\subseteq \N(M\s)$, i.e.,
$\{M\s a^{i_1t}b^2, M\s a^{-i_1t}b^{-2}, Ma^{i_2t}\}\subseteq \N(M\s).$
 Hence, we can conclude that
$$\{ Ma^{x+i_1}b^{y+1},Ma^{x-i_1}b^{y-1},M\s a^{x+i_2}b^y\}\subseteq \N(Ma^xb^y),$$
and
 $$\{ M\s a^{x+i_1t} b^{y+2},M\s a^{x-i_1t}b^{y-2},M a^{x+i_2t}b^y\}\subseteq \N(M\s a^xb^y),$$
for any $x\in\ZZ_p$, $y\in \ZZ_5$.

Assume $i_1\neq 0$. We assert that there exists an integer $\ell$ such that $(a^{i_1}b)^\ell=a^{i_1t}b^2$. If $t=2$, then we can take $\ell=2$ and the asset holds. Now we assume that $t\neq 2$. Since $(p, 5)=1$, there exist integers $m_1$ and $m_2$ such that $pm_1+5m_2=1$ and
$(t-2)(pm_1+5m_2)=t-2$, i.e., $t+p(2-t)m_1=2+5(t-2)m_2$. Taking $n_1=(2-t)m_1i_1$ and $n_2=(t-2)m_2$, one yields that $\ell=t+p(2-t)m_1=2+5(t-2)m_2$ and the assert holds.
Define for each $k'\in\ZZ_{5p}$: $u_{k'}=M(a^{i_1}b)^{k'}$ and $v_{k'}=M\s(a^{i_1}b)^{k'}$.
Then there exists a subgraph $X_1\subseteq X$ such that $V(X_1)=V(X)$ and
$E(X_1)=\{ u_{k'}u_{k'+1}, v_{k'}u_{k'}, v_{k'}v_{k'+\ell}\di k'\in\ZZ_p \}$, where $\ell$ is the integer constructed previously. By construction, $X_1$ is a generalized Petersen graph. Applying Proposition \ref{petersen}, we conclude that $X_1$ admits a Hamilton cycle, and consequently $X$ contains a Hamilton cycle. Thus, Claim 1 holds in this case.

Assume $i_1=0$. Given that $\{Ma^{i_1}b, Ma^{-i_1}b^{-1}, M\s a^{i_2}\}\subseteq\N(M)$ for some $i_1, i_2\in\ZZ_p$ and $B_{5i+j}$ is relabeled by $M\s^i\lg a\rg b^j$ for $i\in\{0, 1\}$ and $j\in\ZZ_5$, we derive the following specific relabelings:
$B_0:=M\lg a\rg$, $B_1:=M\lg a\rg b$, $B_4:=M\lg a\rg b^4$ and $B_5:=M\s\lg a\rg$. Then $Mb\in\N(M)\cap B_1$, $Mb^4\in\N(M)\cap B_4$ and $M\s a^{i_2}\in \N(M)\cap B_5$.
Suppose $|\N(M)\cap B_1|\geq2$ or $|\N(M)\cap B_4|\geq2$. Then there exists a positive integer $\iota$ such that either $Ma^\iota b\in \N(M)$ or $Ma^\iota b^{-1}\in\N(M)$. This situation reduces to the previously analyzed case where $i_1\neq 0$.
Now, consider the complementary case where $|\N(M)\cap B_1|=|\N(M)\cap B_4|=1$, which implies $\N(M)\cap B_1=\{Mb\}$ and $\N(M)\cap B_4=\{Mb^4\}$.

Assume there exists an edge of $X$ contained entirely within $B_0$, that is, $Ma^e\in\N(M)$ for some $e\in\ZZ_p\setminus\{0\}$.
Then $\{Ma^e, Ma^{-e}, Mb, Mb^4\}\subseteq\N(M)$.
Then the graph $Y_0$ contains a subgraph $Y_0'$ defined by the coset construction: 
$Y_0':=\Cos(M\lg a\rg\lg b\rg, M, \cup_{g\in S}MgM)$, where $S=\{a^e,a^{-e},b,b^4\}$. Then $Y_0'$ and $Y_0$ are connected by Proposition \ref{coset-graph-pro}, as $\lg M, a^e,a^{-e},$ $b,b^4 \rg=M\lg a \rg\lg b\rg$.
The connectivity of $X$ guarantees that its quotient graph $X_{\lg ab\rg}$ is connected, which implies the existence of edges between the subgraphs $Y_0$ and $Y_1$. Moreover, since $Y_1=Y_0^\s$, it follows that $Y_1$ is also connected and $d(Y_1)\geqslant4$. By applying Proposition~\ref{deg3}, we conclude that $X$ contains a Hamilton cycle.

\begin{figure}
  \centering
{\tiny \begin{tikzpicture}[node distance=2cm, auto]
    % Nodes
    \node[ellipse, draw] (u0) {$u_{(0,0)}$};
    \node[ellipse, draw, below of=u0] (v0) {$v_{(0,2)}$};
    \node[ellipse, draw, right of=u0] (u1) {$u_{(1,2)}$};
    \node[ellipse, draw, below of=u1] (v1) {$v_{(1,1)}$};
    \node[ellipse, draw, right of=u1] (u2) {$u_{(2,1)}$};
    \node[ellipse, draw, below of=u2] (v2) {$v_{(2,1)}$};
    \node[ellipse, draw, right of=u2] (u3) {$u_{(3,1)}$};
    \node[ellipse, draw, below of=u3] (v3) {$v_{(3,2)}$};
    \node[ellipse, draw, right of=u3] (u4) {$u_{(4,2)}$};
    \node[ellipse, draw, below of=u4] (v4) {$v_{(4,2)}$};
    \node[ellipse, draw, right of=u4] (u5) {$u_{(5,2)}$};
    \node[ellipse, draw, below of=u5] (v5) {$v_{(5,1)}$};
    \node[ellipse, draw, right of=u5] (u6) {$u_{(6,1)}$};
    \node[ellipse, draw, below of=u6] (v6) {$v_{(6,0)}$};

    % Edges
    \draw (u0) -- (v1);
    \draw (v0) -- (u6);
    \draw (v6) -- (u0);
    \draw (u1) -- (v2);
    \draw (u1) -- (v0);
    \draw (u2) -- (v1);
    \draw (u2) -- (v3);
    \draw (u3) -- (v2);
    \draw (u3) -- (v4);
    \draw (u4) -- (v3);
    \draw (u4) -- (v5);
    \draw (u5) -- (v4);
    \draw (u5) -- (v6);
    \draw (u6) -- (v5);
\end{tikzpicture}}
  \caption{Hamilton cycle in $X$ for $|X|=70$ and $X_\B$ is a Petersen graph }\label{c1}
  \end{figure}
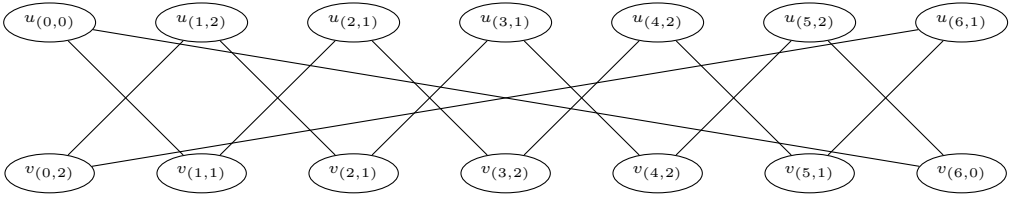

Assume $X$ has no edges within $B_0$, i.e.,  $\N(M)\cap B_0=\emptyset$. Note that $\{Mb, Mb^4, M\s a^{i_2}\}\subseteq\N(M)$, $\N(M)\cap B_1=\{Mb\}$, $\N(M)\cap B_4=\{Mb^4\}$ and $M\s a^{i_2}\in \N(M)\cap B_5$ for some $i_2\in\ZZ_p$. Now rewrite $\ell=r_2$, that is, $M\s a^\ell\in N(M)\cap B_5$. Since $a^\s=a^t$ and $t^{2^{k+1}}\equiv-1 (\bmod\, p)$, $a^{\s^{2^{k+1}}}=a^{t^{2^{k+1}}}=a^{-1}$. From $M\s a^\ell\in\N(M)$ and $M=\lg\s^2\rg$, we derive that $M\s a^\ell\s^{2^{k+1}}\in\N(M\s^{2^{k+1}})=\N(M)$ and
$M\s a^\ell\s^{2^{k+1}}=M(\s^{-2})^{2^k}\s a^\ell\s^{2^{k+1}}=M\s\s^{-2^{k+1}}a^\ell\s^{2^{k+1}}=M\s (a^{\ell})^{\s^{2^{k+1}}}=
M\s a^{-\ell}$. Then $\{Mb, Mb^4, M\s a^\ell, M \s a^{-\ell}\}\subseteq\N(M)$ and so $M\s b^2=M\s b^\s=Mb\s\in\N(M\s)$, $M\s b^3=M\s (b^4)^\s=Mb^4\s\in\N(M\s)$, $M\s\in\N(Ma^{-\ell})$ and $M\s\in\N(Ma^\ell)$. That is to say $\{M\s b^2, M\s b^3, M a^\ell, M a^{-\ell}\}\subseteq\N(M\s)$. Further, we can conclude that
$$\{Ma^xb^{y+1}, Ma^xb^{y+4}, M\s a^{x+\ell}b^y, M\s a^{x-\ell}b^y\}\subseteq N(Ma^xb^y)$$
and
$$\{M\s a^xb^{y+2},M\s a^xb^{y+3}, Ma^{x+\ell}b^y,Ma^{x-\ell}b^y\}\subseteq N(M\s a^xb^y),$$
for any $x\in\ZZ_p$ and $y\in\ZZ_5$.
We define vertex sets with explicit parameterization. Let $u_{(k',\jmath)}=Ma^{\ell k'}b^\jmath$ and $v_{(k',\jmath)}=M\s a^{\ell k'}b^\jmath$, where $k'\in\ZZ_p$ and $\jmath\in\ZZ_5$. Then $V(X)=\{u_{(k',\jmath)},v_{(k',\jmath)}\mid k'\in\ZZ_p, \jmath\in\ZZ_5\}$.
Let $A=\{ u_{(k',0)}, u_{(k',1)}, u_{(k',2)}, u_{(k',3)}, u_{(k',4)}\}$ and $B=\{v_{(k'+1, 0)}, v_{(k'+1, 1)}, v_{(k'+1, 2)},$ $v_{(k'+1, 3)}, v_{(k'+1, 4)}\}$.
Then the neighborhood of $Ma^{x}b^y$ and $M\sigma a^xb^y$ for $x\in\ZZ_p$ and $y\in\ZZ_5$ leads to that $X[A,B]$ has a subgraph which is isomorphic to Petersen graph.
We set
$$W_{(k',0,1)}:=u_{(k',0)}u_{(k',1)}u_{(k',2)}u_{(k',3)}u_{(k',4)}
v_{(k'+1,4)}v_{(k'+1,2)}v_{(k'+1,0)}v_{(k'+1,3)}v_{(k'+1,1)},$$
$$W_{(k',1,2)}:=u_{(k',1)}u_{(k',2)}u_{(k',3)}u_{(k',4)}u_{(k',0)}
v_{(k'+1,0)}v_{(k'+1,3)}v_{(k'+1,1)}v_{(k'+1,4)} v_{(k'+1,2)},$$
$$W_{(k',2,1)}:=u_{(k',2)}u_{(k',1)}u_{(k',0)}u_{(k',4)}u_{(k',3)}
v_{(k'+1,3)}v_{(k'+1,0)}v_{(k'+1,2)}v_{(k'+1,4)}v_{(k'+1,1)},$$
$$W_{(k',1,0)}:=u_{(k',1)}u_{(k',0)}u_{(k',4)}u_{(k',3)}u_{(k',2)}
v_{(k'+1,2)}v_{(k'+1,4)}v_{(k'+1,1)}v_{(k'+1,3)}v_{(k'+1,0)}.$$
Then
$$\begin{array}{ll}
W_{(0,0,1)}\sim W_{(2,1,2)}\sim W_{(4,2,1)}\sim\cdots\sim
W_{(p-4,1,2)}
\sim W_{(p-2,2,0)}\sim W_{(0,0,1)}
\end{array}$$
is a Hamilton cycle (see Figure~\ref{c1} for example), as $v_{(p-1,0)}\sim u_{(0,0)}$. Thus, Claim 1 holds.

\vskip 3mm
{\it Case 2: $\C_G(P)\cong\ZZ_{5p}:\ZZ_2$ or $\ZZ_{5p}:\ZZ_4$}
\vskip 3mm
By Lemma \ref{Z5Z4}, we conclude that $G=\lg a\rg\times\lg b\rg:\lg\s\rg=\lg a,b,\sigma\mid a^p=b^5=\sigma^4=1, [a,b]=1, a^\sigma\in\lg a\rg, b^\sigma=b^2\rg\cong\ZZ_{5p}:\ZZ_4$ and $K=\lg a \rg$. Note that $\C_G(\lg a\rg)\cong\ZZ_{5p}:\ZZ_2$ or $\ZZ_{5p}:\ZZ_4$, we conclude that
$$G=\lg a,b,\s\di a^p=b^5=\s^4=1,[a,b]=1,a^\s=a^\iota,b^\s=b^2\rg,$$
where $\iota=-1$ or 1 depends on $C_G(\lg a\rg)\cong\ZZ_{5p}:\ZZ_2$ or $\ZZ_{5p}\rtimes\ZZ_4$, respectively.
Noting that $\N(B_0)=\{ B_1,B_4,B_5\}$ (as shown in Figure \ref{coPetersengraph}) and the relabeling: $B_{5i+j}:=M\s^i\lg a\rg b^j$ for $i\in\{0, 1\}$ and $j\in\ZZ_5$, we conclude that $Ma^\ell b\in\N(M)$ for some $\ell\in \ZZ_p$.

Assume $\ell\neq0$. Followed from $B_5\in \N(B_0)$ that we can take $M\s a^\jmath\in\N(M)$ for some $\jmath\in\ZZ_p$.
Then $\{Ma^\ell b, Ma^{-\ell} b^{-1}, M\s a^\jmath\}\subseteq\N(M)$ and $\{ Ma^\ell b\s^2, Ma^{-\ell} b^{-1}\s^2\}\subseteq\N(M\s^2)=\N(M)$ with $Ma^\ell b\s^2=M\s^2\s^{-2}a^\ell \s^2\s^{-2}b\s^2=M((a^\s)^\s)^\ell (b^\s)^\s=Ma^\ell b^{-1}$ and $Ma^{-\ell} b^{-1}\s^2=Ma^{-\ell}b$. It follows that $\{Ma^\ell b, Ma^{-\ell} b^{-1}, M\s a^\jmath, Ma^\ell b^{-1}, Ma^{-\ell}b\}\subseteq\N(M)$
Then the graph $Y_0$ has a subgraph $Y_0':=\Cos(M\lg a\rg\lg b\rg, M, \cup_{g\in S_1} MgM)$, where $S_1=\{a^\ell b,a^{-\ell} b^{-1},a^\ell b^{-1},a^{-\ell}b\}$. Then $Y_0'$ and $Y_0$ are connected by Proposition \ref{coset-graph-pro}, as $\lg M, a^\ell b,a^{-\ell} b^{-1},a^\ell b^{-1},a^{-\ell}b\rg=M\lg a \rg\lg b\rg$.
And $Y_0$ is adjacent to $Y_1$, as $M\s a^\jmath\in\N(M)$. Since $Y_1=Y_0^\s$, it follows that $Y_1$ is also connected and $d(Y_1)\geqslant4$. By Proposition~\ref{deg3}, we conclude that $X$ contains a Hamilton cycle.

Assume $\ell=0$. With the same arguments as case $i_1=0$ in Case 1, we can only consider the case $\N(M)\cap B_1=\{Mb\}$ and $\N(M)\cap B_4=\{Mb^4\}$. Suppose that $Ma^{j_1}\in\N(M)$ for some $j_1\in\ZZ_p$. Then $\{Ma^{j_1}, Ma^{-j_1}, Mb, Mb^4\}\subseteq\N(M)$.
Then $\Cos(M\lg a\rg\lg b\rg,M,$ $\cup_{g\in S_2} MgM)$ is a connected subgraph of $Y_0$ with degree greater than $4$,
where $S_2=\{a^{j_1},a^{-j_1},$ $b,b^{-1}\}$.
 Applying Proposition $\ref{deg3}$, we conclude that $X$ contains a Hamilton cycle.
So in what follows, we assume $\N(M)\cap (\bigcup_{j=0}^4B_j)=\{Mb,Mb^4\}$.
Derived from the connectivity of graph $X$, there exist distinct $l_1,l_2\in\ZZ_p$ with $l_1\neq 0$ such that $M\s a^{l_1}\in \N(M)\cap B_5$ if $\C_G(\lg a\rg)\cong\ZZ_{5p}\rtimes\ZZ_4$, and $\{M\s a^{l_1}, M\s a^{l_2}\}\subseteq \N(M)\cap B_5$ if $\C_G(\lg a\rg)\cong\ZZ_{5p}\rtimes\ZZ_2$.
Then we get
$$\{M\s^ia^xb^{y+1+i},M\s^ia^xb^{y+4-i},M\s^{1-i} a^{x+l_1}b^y,M\s^{1-i} a^{x-l_1}b^y\}\subseteq \N(M\s^ia^xb^y)$$
if $\C_G(\lg a\rg)\cong\ZZ_{5p}\rtimes\ZZ_4$, and
$$\{M\s^ia^xb^{y+1+i},M\s^ia^xb^{y+4-i},M\s^{1-i} a^{x+(-1)^il_1}b^y,M\s^{1-i} a^{x+(-1)^il_2}b^y\}\subseteq \N(M\s^ia^xb^y)$$ if $\C_G(\lg a\rg)\cong\ZZ_{5p}\rtimes\ZZ_2$, for $i=0,1$ and $x\in\ZZ_p$, $y\in\ZZ_5$. By an argument entirely analogous to Case 1, we conclude that $X$ contains a Hamilton cycle.
\vskip 3mm
\textbf{Claim 2:} If $X_\B$ is isomorphic to the complement of the Petersen graph, then $X$ contains a Hamilton cycle.
\vskip 3mm
Suppose $X_\B$ is isomorphic to the complement of the Petersen graph. Then one can get a subgraph of $X_\B$ which is isomorphic to the Petersen graph say $X_1$ (see $\Gamma_2$ in Figure~\ref{coPetersengraph}).
By an argument completely analogous to the proof of Claim 1, we conclude that $X$ contains a Hamilton cycle. Therefore, Claim 2 holds.
\vskip 3mm
\textbf{Claim 3:} If $X_\B\cong\K_{5,5}-5\K_{1,1}$, then $X$ contains a Hamilton cycle.
\vskip 3mm
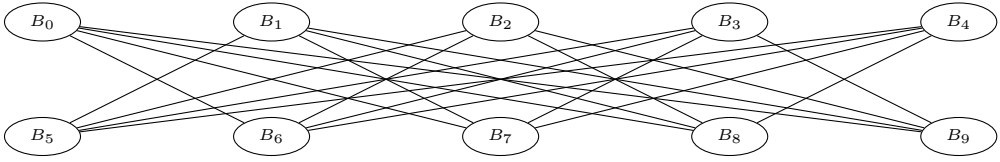
\begin{figure}\centering
{\tiny\begin{tikzpicture}[
  node distance=1cm and 2cm,
  every node/.style={draw, ellipse, minimum width=1cm, minimum height=0.5cm}
]
% Nodes
\node (A0) {$B_0$};
\node (A1) [right=of A0] {$B_1$};
\node (A2) [right=of A1] {$B_2$};
\node (A3) [right=of A2] {$B_3$};
\node (A4) [right=of A3] {$B_4$};
\node (B0) [below=of A0] {$B_5$};
\node (B1) [right=of B0] {$B_6$};
\node (B2) [right=of B1] {$B_7$};
\node (B3) [right=of B2] {$B_8$};
\node (B4) [right=of B3] {$B_9$};
% Edges
\foreach \i in {0,1,2,3,4} {
  \foreach \j in {0,1,2,3,4} {
    \ifnum\i=\j
      \else
        \draw (A\i) -- (B\j);
    \fi
  }
}
\end{tikzpicture}}
  \caption{$K_{5,5}-5K_{1,1}$}\label{K55-K11}
  \end{figure}

Assume $X_\B\cong\K_{5, 5}-5\K_{1, 1}$  (see Figure \ref{K55-K11} for example). Then $X_\B$ is $G/K$ edge-transitive.
From Figure~\ref{K55-K11}, we observe that $\N(B_0)=\{ B_6,B_7,B_8,B_9\}$ in $X_\B$.
Further, assuming $d(B_0,B_6)\geqslant2$, the edge-transitivity of $X_\B$ implies that $d(B_{j'},B_{5+j})\geqslant2$ for distinct $j',j\in\ZZ_5$. Since $a$ is a $(10,p)$-semiregular automorphism, Proposition~\ref{pro:4} directly implies that $X$ contains a Hamilton cycle, as desired.
For subsequent arguments, we now assume the complementary case: $d(B_0, B_6)=1$. Then $d(B_{j'}, B_{5+j})=1$ for distinct $j', j\in\ZZ_5$. Furthermore, the connectivity of graph $X$ guarantees $\N(M)\cap B_0\neq\emptyset$.
From the group structure of $G$ in Lemma \ref{Z5Z4}, we recall that $|\sigma|=2^{k+2}$.

Assume $k\neq0$, then $\C_G(P)=\lg ab\rg\cong\ZZ_{5p}$, $G=\lg a,b,\s\di a^p=b^5=\s^{2^{k+2}}=1,[a,b]=1,a^\s=a^t,b^\s=b^2, t^{2^{k+1}}\equiv-1 (\bmod\, p)\rg$ and $K=\lg a,c\mid c=\sigma^4, a^c=a^{t^4}, t^{2^{k+1}}\equiv-1 (\bmod\, p) \rg$. In particular, $a^{c}=a^{\sigma^4}\neq a$, as $\C_G(P)\cap K=\lg a\rg$ and then $a^{\s^2}\neq a^{-1}$.
Set $Ma^\jmath\in\N(M)\cap B_0$ with $\jmath\in\ZZ_p\setminus\{0\}$. Noting that $Ma^\jmath\sigma^2\in\N(M\sigma^2)=\N(M)$, $B_0=M\lg a\rg$, $\lg a\rg^{\sigma^2}=\lg a\rg$ and $a^{\s^2}\neq a^{-1}$, one yields that $\{Ma^\jmath, Ma^{-\jmath}, Ma^\jmath\s^2, Ma^{-\jmath}\s^2\}\subseteq \N(M)\cap B_0$. Further, we derive from the subgraph $X\lg B_0\rg$ is connected and  $d( B_0)\geqslant4$ and that $X_\B$ admits a Hamilton cycle. Applying Proposition~\ref{deg3}, $X$ contains a Hamilton cycle, as desired.

Throughout the remainder of this proof, we fix $k=0$. By the arguments of Case 2 of Claim 1, we deduce that $G=\lg a,b,\s\di a^p=b^5=\s^4=1,[a,b]=1,a^\s=a^\iota,b^\s=b^2\rg,$
where $\iota=-1$ or 1 depends on $C_G(\lg a\rg)\cong\ZZ_{5p}:\ZZ_2$ or $\ZZ_{5p}\rtimes\ZZ_4$, respectively.
 Take $M\s a^\ell b\in\N(M)\cap B_6$ for some $\ell\in\ZZ_p$.
Then $M\s a^{-\iota\ell}b^2=M(\s a^{\ell}b)^{-1}\in\N(M)$, $M\s a^{\ell}b^{4}=M\s a^\ell b\sigma^2\in\N(M\sigma^2)=\N(M)$ and $M\s a^{-\iota\ell}b^3=M(\s a^{\ell}b)^{-1}\s^2\in\N(M\sigma^2)=\N(M)$, that is $\{M\s a^{\ell}b,M\s a^{\ell}b^{4},M\s a^{-\iota\ell}b^2,M\s a^{-\iota\ell}b^3\}\subseteq\N(M).$ And $\{M a^{-\ell}b,M a^{\iota\ell}b^{2},M a^{\iota\ell}b^3,M a^{-\ell}b^4\}\subseteq\N(M\sigma).$
So, for $i\in\{0,1\}$, $x\in\ZZ_p$ and $y\in\ZZ_5$, we can easily arrive at the following conclusion:
$$\{M\s^{1-i} a^{(-1)^i\ell+x}b^{1+y},M\s^{1-i} a^{(-1)^i(-\iota\ell)+x}b^{2+y},M\s^{1-i} a^{(-1)^i(-\iota\ell)+x}b^{3+y},$$
$$M\s^{1-i} a^{(-1)^i\ell+x}b^{4+y}\}\subseteq \N(M\s^i a^xb^y).$$
Since $\N(M)\cap B_0\neq \emptyset$ and $B_0$ can be relabeled as $M\lg a \rg$, we may choose $Ma^{\ell'}\in \N(M)\cap B_0$ for some $\ell'\in\ZZ_p\setminus\{0\}$.
Then $\{ Ma^{\ell'}, Ma^{-\ell'}\}\subseteq \N(M)$.
Then we can also get $\{M\sigma^i a^{\ell'+x}b^y,M\sigma^i a^{-\ell'+x}b^y\}\subseteq \N(M\sigma^ia^xb^y)$ for $i=0,1$ and any $x\in\ZZ_p$, $y\in\ZZ_5$.
For each vertex induced graph $X\lg B_{5i+j} \rg$, there exists a Hamilton path
{\small$$\begin{array}{lll}
\mathcal{C}_{5i+j}(M\sigma^ia^{x}b^j,M\sigma^ia^{x+(p-1)\ell'}b^j)&:=&M\sigma^ia^xb^j\sim M\sigma^ia^{\ell'+x}b^j\sim M\sigma^ia^{2\ell'+x}b^j\sim \\
&&\cdots\sim M\sigma^ia^{(p-1)\ell'+x}b^j,
\end{array}$$}
and define
$${\small\begin{array}{lll}
\mathcal{C}_{5i+j}^{-1}(M\sigma^ia^{x}b^j,M\sigma^ia^{x+(p-1)\ell'}b^j)&:=&M\sigma^ia^{(p-1)\ell'+x}b^j\sim \cdots\sim M\sigma^ia^{2\ell'+x}b^j\\
&&\sim M\sigma^ia^{\ell'+x}b^j\sim M\sigma^ia^xb^j,
\end{array}}$$
\f where $i=0,1$ and for any $x\in\ZZ_p$ and $j\in\ZZ_5$.
Then we construct a Hamilton cycle in $X$ as follows:
 {\small$$
\begin{array}{ll}
\mathcal{C}_{0}(M,Ma^{-\ell'})\sim \mathcal{C}_{5+1}^{-1}(M\sigma a^{\ell}b ,M\sigma a^{\ell-\ell'}b )\sim\mathcal{C}_{0+2}(Mb^2,Ma^{-\ell'}b^2) \sim\mathcal{C}_{5+3}^{-1}(M\sigma a^{\ell}b^3 ,M\sigma a^{\ell-\ell'}b^3)\\ \sim\mathcal{C}_{0+4}(Mb^4,Ma^{-\ell'}b^4)\sim
\mathcal{C}_{5+0}^{-1}(M\sigma a^{\ell} ,M\sigma a^{\ell-\ell'})\sim \mathcal{C}_{0+1}(Mb,Ma^{-\ell'}b)
\sim\\
\mathcal{C}_{5+2}^{-1}(M\sigma a^{\ell}b^2 ,M\sigma a^{\ell-\ell'}b^2)
\sim\mathcal{C}_{0+3}(Mb^3,Ma^{-\ell'}b^3)\sim
\mathcal{C}_{5+4}^{-1}(M\sigma a^{\ell}b^4 ,M\sigma a^{\ell-\ell'}b^4),
\end{array}$$}
as $M\sigma a^\ell b^4\in \N(M)$.
\qed

\begin{lemma}\label{A5}
Using the notations of Theorem \ref{p-blocks}, if $G/K\cong\A_5$, then the graph $X$ contains a Hamilton cycle.
\end{lemma}
\demo Assume that $G/K\cong\A_5$. Since $K\cong\ZZ_p:\ZZ_l$ where $l\di p-1$, we have $G\cong (\ZZ_p:\ZZ_l).\A_5$. We derive from $G/\C_G(\lg a\rg)$ is cyclic, $\C_G(\lg a\rg)\cap K=\lg a\rg$ and $(\C_G(\lg a\rg)K)/K\trianglelefteq G/K$ that  $\C_G(\lg a\rg)=\lg a\rg\times T$ with $T\cong\A_5$. Since $\lg a\rg$ acts transitively on each $K$-block and $T$ acts transitively on the set of $K$-blocks, it follows that $\C_G(\lg a\rg)$ acts transitively on $V(X)$. That is $G$ contains a transitive subgroup $\C_G(\lg a\rg)\cong\ZZ_p\times \A_5$ acting on $V(X)$. From the minimal transitivity of $G$, we can draw that $G=\C_G(\lg a\rg)$ and $K=\lg a\rg$. On the other hand, since $T\cong\A_5$ acts on $10$ popints with three suborbits of length $1$, $3$ and $6$, it follows that $X_\B$ must be either the Petersen graph or the complement of the Petersen graph, unless $X_{\B}\cong K_{10}$. Thus, we may restrict our analysis to the cases where $X_\B$ is the Petersen graph, or its complement.
Observe that the complement of the Petersen graph contains a subgraph isomorphic to the Petersen graph itself (see $\Gamma_2$ in Figure \ref{coPetersengraph}).Therefore, without loss of generality, we may assume $X_{\B}$ is the Petersen graph in our subsequent analysis.

Pick up $b\in T$ such that $o(b)=5$. Since $o(a)=p$ (where $p\geqslant7$ is prime) and $[a, b]=1$, the subgroup $\ZZ_{5p}\cong\lg ab\rg$. Moreover, $ab$ is a semiregular automorphism of graph $X$, with exactly two orbits $\Delta_0$ and $\Delta_1$ on $V(X)$. For $\jmath\in\{0,1\}$, let $Y_\jmath=X\lg\Delta_\jmath\rg$ be the induced subgraph. Then $V(Y_\jmath)=\{M\s^\jmath a^xb^y\di x,y\in\ZZ\}$, where $\jmath=0,1$, $\s\in T\setminus M\lg a b\rg$ and $M$ is a point stabilizer of $G$ on $V(X)$ such that $M\cong\S_3$ with $\sigma^2\in M$.
Now we relabel $B_{5\jmath+j}$ as $M\s^\jmath\lg a\rg b^j$ for $\jmath\in\{0, 1\}$ and $j\in\ZZ_5$. Then the graph $X_\B$ is the  $\Gamma_1$ in the Figure~\ref{coPetersengraph} (up to isomorphic). Further, a detailed analysis of graph $\Gamma_1$ in Figure~\ref{coPetersengraph} reveals that $\N(B_0)=\{ B_1,B_4,B_5\}$ in $X_\B$. Consequently, we may select an element $Ma^\ell b\in\N(M)\cap B_{1}$ for some $\ell\in\ZZ_p$. This implies that $Ma^{-\ell}b^{-1}=M(a^\ell b)^{-1}\in\N(M)$.
Consider the right multiplication of $M$ on $[T:M]$. Since $M\cong\S_3$ and $o(b)=5$, one yields that $Mb\neq Mb^{-1}$, in particular $M\sigma\not\in\{Mb, Mb^{-1}\}$. Then one yields that $MbM=Mb\cup Mb^{-1}\cup M\s$ and so there exit elements $\tau_1$ and $\tau_2$ in $M$ such that $Mb\tau_1= Mb^{-1}$ and $Mb\tau_2= M\s$. Therefore, $Ma^\ell b^{-1}=Ma^\ell b\tau_1\in\N(M\tau_1)=\N(M)$ and $Ma^{-\ell}b=M(a^\ell b^{-1})^{-1}\in\N(M)$.

Assume $\ell\neq0$. Then
$\{ Ma^{\ell}b, Ma^{\ell}b^{-1}, Ma^{-\ell}b^{-1}, Ma^{-\ell}b\}\subseteq \N(M)$. Furthermore, this implies that
$\{ Ma^{x+\ell}b^{y+1},Ma^{x+\ell}b^{y-1},Ma^{x-\ell}b^{y-1},Ma^{x-\ell}b^{y+1}\}\subseteq \N(Ma^xb^y)$ for $x\in\ZZ_p$ and $y\in\ZZ_5$.
Since $\lg a^\ell b\rg\cong\ZZ_{5p}$, one yields that $Y_0$ is connected and $\deg(Y_0)\geqslant4$.
Applying Proposition~\ref{deg3}, we deduce that $X$ contains a Hamilton cycle, as desired.

Assume $\ell=0$. Then $Mb\in\N(M)$ and $M\s=Mb\tau_2\in\N(M\tau_2)=\N(M)$. Consequently, we have that $\{Mb, Mb^{-1}, M\sigma, M\sigma^{-1}\}\subseteq\N(M)$.
Observe that $X$ can be represented as a coset graph $\Cos(G,M,\cup_{g\in S}MgM)$, where $S\subseteq G$ satisfies $\cup_{g\in S}Mg^{-1}M=\cup_{g\in S}MgM$, and $\lg M,S \rg=G$. 
Given that $\{Mb, Mb^{-1}, M\sigma, M\sigma^{-1}\}\subseteq\N(M)$, it follows that $M\{b,b^{-1},\sigma,\sigma^{-1}\}M\subset MSM$ and so $M\sigma^\jmath a^{\ell'} b^s\in\N(M)$ for some $\jmath\in\{0,1\}$, $\ell'\neq 0$ and $s\in\ZZ_5$, as $\lg M,S\rg=G$. 
If $\jmath=1$, then  $M\sigma a^{\ell'} b^s\in\N(M)$ for some $\ell'\neq 0$ and $s\in\ZZ_5$. By computing the suborbit of $M$ acting on the point $M\sigma a^{\ell'} b^s$, we obtain that $M\sigma a^{\ell'} b^sM=M\sigma M b^{s'} a^{\ell'}=Mb M b^{s'} a^{\ell'}$, for some integer $s'$, as $MbM=M\sigma M$ and $\lg b \rg M=M\lg b\rg$. Now, we get that $M\sigma a^{\ell'} b^sM=M b^{s'}a^{\ell'}M$ for some integer $s'$. Then either $M a^{\ell'}b^{\pm 1}\in \N(M)$ or $M a^{\ell'}\in \N(M)$.

If either $M a^{\ell'}b^{\pm 1}\in \N(M)$ or $\jmath=0$, then the existence of a Hamilton cycle in $X$ follows from the $\ell\neq 0$ case previously established. Now, we only need to consider the case that $\{Mb, Mb^{-1}, Ma^{\ell'}, Ma^{-\ell'}\}\subseteq\N(M)$ with $\ell'\neq 0$.
Then $d(Y_0)\geqslant4$ and $Y_0$ is connected as $\lg a^\ell b\rg\cong\ZZ_{5p}$. 
By Proposition~\ref{deg3} again, one yields that $X$ contains a Hamilton cycle, as desired.
\qed

\begin{lemma}\label{Z10gp}
Using the notations of Theorem \ref{p-blocks},
if $G\cong K.\D_{10}$, then $G\cong\ZZ_{5p}:\ZZ_{2l}$ where $l$ is $2$-power;
if $G\cong K.\ZZ_{10}$, then $G\cong\ZZ_p:\ZZ_{10l}$.
\end{lemma}
\demo Note that $K=\lg a\rg\rtimes\lg\s\rg\cong\ZZ_p:\ZZ_l$ is of affine type and $K\leq\AGL(1, p)$.
Then $P=\lg a\rg\cong\ZZ_p$ is a Sylow $p$-subgroup of $G$, where $p\geqslant7$,
and $P\unlhd G$.

First, consider the case $G\cong K.\D_{10}$.
Since $G/K\cong \D_{10}$ and $G/\C_G(P)$ is cyclic, the group $G/\C_G(P)K$ is a quotient of $\D_{10}$ and is cyclic,
which implies that either $G=\C_G(P)K$ or $|G:\C_G(P)K|=2$.
Suppose that $G=C_G(P)K$. Then $G=\C_G(P)K=\C_G(P)\lg a\rg\lg \s\rg=\C_G(P)\lg \s\rg$ and
$10pl=|G|=|\C_G(P)\lg \s\rg|=|\C_G(P)||\lg \s\rg|$, as $\C_G(P)\cap\lg \s\rg=1$.
Since $\lg\s\rg$ is a point stabilizer, by the minimality of $G$,
we get that $G=\C_G(P)$, $K=K\cap \C_G(P)=\lg a\rg$ and $|G|=|\C_G(P)|=10p$.
By Proposition~\ref{complement2}, there exists a subgroup $I$ of $G$ such that
$G=\lg a\rg I$, $I\cap\lg a\rg=1$ and $I\cong \D_{10}$.
Then $G=\lg a\rg\times I\cong\ZZ_p\times \D_{10}\cong \ZZ_{5p}\rtimes\ZZ_2$, as desired.
Suppose that $|G:\C_G(P)K|=2$.
Then $5pl=|\C_G(P)K|=|\C_G(P)\lg \s\rg|=|\C_G(P)||\lg \s\rg|$, which implies $\C_G(P)\cong\ZZ_{5p}$.
Let $c$ be a $2$-element in $G\setminus \C_G(P)K$.
Then $\C_G(P)\cap\lg c\rg=1$. Set $L=\C_G(P)\lg c\rg$ and $\C_G(P)K<LK\leq G$.
By $|G:\C_G(P)K|=2$ and $LK=L\lg\s\rg$, we get $G=LK=L\lg\s\rg$.
By the minimality of $G$ again, we get $G=L=\C_G(P)\lg c\rg=\C_G(P)\rtimes\lg c\rg$,
which implies $|G|=|\C_G(P)\lg c\rg|=5p\cdot\o(c)$.
Then $\o(c)=2l$, where $l$ is $2$-power, and so $G\cong\ZZ_{5p}:\ZZ_{2l}$.

Second, consider the case $G\cong(\ZZ_p:\ZZ_l).\ZZ_{10}$.
By Proposition~\ref{complement2}, there exists a subgroup $I$ of $G$
such that $G=\lg a\rg I$, $\s\in I$ and $I\cap \lg a\rg=1$.
Since $I\cong \lg a\rg I/\lg a\rg=G/\lg a\rg\cong \ZZ_l.\ZZ_{10}$, we set
$I=\lg \s,b\di \s^l=1,b^{10}=\s^\imath\rg=\lg b\rg \lg\s\rg$, for some integer $\imath$.
Note that $\lg\s\rg$ is a point stabilizer of $V(X)$ and
$G=\lg a\rg I=\lg a\rg (\lg b\rg \lg\s\rg)=(\lg a\rg \lg b\rg) \lg\s\rg=\lg a,b\rg\lg\s\rg.$
Then $\lg a,b\rg=\lg a\rg\rtimes\lg b\rg$ acts transitively on $V(X)$.
By the minimality of $G$, we get $l=\o(\s^\imath)$ and $G=\lg a,b\rg=\lg a\rg\rtimes\lg b\rg\cong\ZZ_p\rtimes\ZZ_{10l}$.
\qed

\begin{lemma}\label{D10}
 Using the notations of Theorem \ref{p-blocks}, if $G/K\cong \D_{10}$, then $X$ contains a Hamilton cycle.
\end{lemma}
\demo
Assume $G/K\cong \D_{10}$. Note that $\ZZ_p:\ZZ_l\cong K\leqslant\AGL(1, p)$, $G\cong(\ZZ_p:\ZZ_l).\D_{10}$ and $G/\ZZ_p\cong\ZZ_l.\D_{10}$ is transitive on $10$ points.
By Lemma \ref{Z10gp}, $G\cong\ZZ_{5p}\rtimes\ZZ_{2l}$, where $l$ is $2$-power.
If $l=1$, then $|G|=10p=|V(X)|$, i.e., $X$ is a Cayley graph of $G$, and so we derive from Proposition~\ref{10pcayley} that $X$ contains a Hamilton cycle, as desired.
For the remainder of the argument, we assume $l>1$. Specifically, we set $l=2^s$ where $s$ is an integer with $s\ge1$. Then
$G=\lg a,b, c\di a^p=b^5= c^{2^{s+1}}=1,a^b=a,b^ c=b^{-1},a^{ c}=a^t,
t^{2^s}\equiv-1\pmod p, 2^{s+1}\di p-1\rg$.
Set $\mathcal{B}$ be a complete $K$-block system of $G$ on $V(X)$ and $M=\lg  c^2\rg$, the point stabilizer. Then
$\B=\{B_{5i+j}:=M c^i b^j\lg a\rg\di i\in\ZZ_2,j\in\ZZ_5\}$.
Then the quotient graph $X_\B$ is a Cayley graph of $G/K$ and so $X_\B$ is not a Petersen graph.
By Proposition \ref{pq-H-cycle}, one yields that  $X_\B$ contains a Hamilton cycle.

Let $Mc^\imath b^\jmath a^k\in \N(M)$, for some $\imath\in\ZZ_2,~\jmath\in\ZZ_5$ and $k\in\ZZ_p$.
Suppose that $(\imath,\jmath)\neq(0,0)$ and $k\neq0$.
Note that $X_{\B}$ is a Cayley graph of $\D_{10}$.
By analyzing, when we consider the action of $\D_{10}$ on the edge set of $X_{\B}$, there exist at least two orbits $\mathcal{E}_A,\mathcal{E}_B$ of the edge set, where $\mathcal{E}_A$ consists of two cycle with vertex set $\{B_{\jmath}\mid \jmath\in\ZZ_5\}$ and $\{B_{5+\jmath}\mid \jmath\in\ZZ_5\}$, and the edge induced subgraph $X_{\B}\lg \mathcal{E}_B\rg\cong 5K_2$. 
Then it is easy to check that there exists a Hamilton cycle of $X_\B$ containing edge $\{B_0,B_{5\imath+\jmath}\}$.
Then we get that $\deg(B_0,B_{5\imath+\jmath})\geqslant2$ as $Mc^{\imath}b^{\jmath}a^k\neq Mc^{\imath}b^{\jmath}a^{kt^2}=Mc^{\imath}b^{\jmath}a^{k}c^2\in \N(M)$, and further, by Proposition~\ref{pro:4}, one yields that $X$ admits a Hamilton cycle, as desired.
So in what follows, we assume that $Mc^\imath b^\jmath a^k\in \N(M)$, only occurs for either $(\imath,\jmath)=(0,0)$ or $k=0$.
According to the connectivity of $X$, we may set $Ma^{k'}\in \N(M)$ and $Mc^{\imath'}b^{\jmath'}$ with some integer $k'\neq 0$ and $(\imath,\jmath)\neq(\imath',\jmath')$.
Then we can get a Hamilton cycle $\mathcal{A}_1$ in the graph $X\lg M\lg c, b\rg \rg :=X\lg \{Mc^ib^j\mid i\in\ZZ_2,j\in\ZZ_5\} \rg$.
And we also get a Hamilton cycle of $X\lg B_0\rg$: $\mathcal{A}_2=MMa^{k'}Ma^{2k'}\cdots Ma^{(p-1)k'}M$,
which implies that $d(B_0)\ge2$ and $B_0$ is connected.
If $d(B_0)\geq 3$, then by Proposition~\ref{deg3}, $X$ admits a Hamilton cycle.
Suppose that $d(B_0)= 2$, that is $\N(M)\cap B_0=\{Ma^{k'}, Ma^{-k'}\}$.
Since $Ma^{k'}c^2=Ma^{k't^2}\in \N(M)\cap B_0$ , we get that either $Ma^{k't^2}=Ma^{k'}$ or $Ma^{k't^2}=Ma^{-k'}$.
For $Ma^{k't^2}=Ma^{k'}$, we get $a^{k'(t^2-1)}\in M$, which implies $t^2-1\equiv 0\pmod p$ as $M\cap\lg a\rg=1$, a contradiction.
For $Ma^{k't^2}=Ma^{-k'}$, we get $a^{k'(t^2+1)}\in M$, which implies $t^2+1\equiv 0\pmod p$.
Then $s=1$, which implies that
$$G=\lg a,b, c\di a^p=b^5= c^4=1,a^b=a,b^ c=b^{-1},a^ c=a^t,t^2\equiv-1\pmod p\rg,$$ and
$\{Mc^ib^ya^{k't^i+x}, Mc^ib^ya^{-k't^i+x}\}\subseteq \N(Mc^ib^ya^x),$ for any $i\in\{0,1\}$ $x\in\ZZ_p$ and $y\in\ZZ_5$.

Set $\mathcal{C}=\{C_{i,j'} \di C_{i,j'}=M c^ia^{j'}\lg b\rg, i\in\ZZ_2, j'\in\ZZ_p\}$, the orbits of $\lg b\rg$.
Consider the quotient graph $X_{\mathcal{C}}$ of order $2p$.
Since $X\lg C_{0,0}\cup C_{1,0}\rg=X\lg M\lg b,c\rg\rg$ contains the cycle $\mathcal{A}_1$,
we get $C_{1,0}\in \N(C_{0,0})$.
Then $\N(C_{0,0})=\{C_{1,0},C_{0,k'},C_{0,-k'}\}$.
Since the subgraph of $X_{\mathcal{C}}$ whose edge set is $\{\{C_{0,0},C_{0,k'}\},\{C_{0,0},C_{0,-k'}\}\}^G$ is not connected and $X_{\mathcal{C}}$ contains a Hamilton cycle, we get that there exists a Hamilton cycle of $X_{\mathcal{C}}$ containing the edge $\{C_{0,0},C_{1,0}\}$.
If $d(C_{0,0},C_{1,0})\ge2$, then by Proposition~\ref{pro:4}, $X$ contains a Hamilton cycle.
So in what follows, we assume that $\deg(C_{0,0},C_{1,0})=1$.
Since $X\lg C_{0,0}\cup C_{1,0}\rg$ contains a Hamilton cycle $\mathcal{A}_1$,
we get that $d( C_{0,0})\ge2$ and $C_{0,0}$ is connected.
If $d( C_{0,0})\ge3$, then $X$ admits a Hamilton cycle by Proposition~\ref{deg3}.
If $d( C_{0,0})=2$, then set $Y_i= M c^i\lg ab\rg$, where $i=0,1$.
Then $d( Y_0)\geq4$ and $X\lg Y_0\rg$ is connected.
By Proposition \ref{deg3} again, $X$ contains a Hamilton cycle, as desired.
\qed

\begin{lemma}\label{Z10}
 Using the notations of Theorem \ref{p-blocks}, if $G/K\cong \ZZ_{10}$, then $X$ contains a Hamilton cycle.
\end{lemma}
\demo
Assume $G/K\cong\ZZ_{10}$. By Lemma \ref{Z10gp}, one yields that $G=\lg a\rg:\lg b\rg\cong \ZZ_p:\ZZ_{10l}$ and $K=\lg a \rg:\lg \sigma \rg\cong \ZZ_p:\ZZ_l$, where $b^{10}=\sigma^\imath$ is of order $l$, for some integer $\imath$. It is clear that $\lg b^{10}\rg\cap \C_G(\lg a\rg)=1$. If $l=1$, then $|G|=10p=|V(X)|$, i.e., $X$ is a Cayley graph of $G$. So we derive from Proposition~\ref{10pcayley} that $X$ contains a Hamilton cycle, as desired.
For the remainder of the argument, we assume $l>1$.
Then $G=\lg a,b\di a^p=b^{10l}=1, a^b=a^t, t^{10}\not\equiv 1\pmod p\rg$.

Let $M=\lg b^{10}\rg$ be a point stabilizer of the action of $G$ on $V(X)$, and $\B=\{B_i\di B_{i}=Mb^i\lg a\rg, i\in\ZZ_{10}\}$, a complete $K$-block system of $G$.
Since $X_\B$ is a Cayley graph of $G/K$ and consequently not isomorphic to the Petersen graph. By Proposition \ref{pq-H-cycle}, one yields that  $X_\B$ contains a Hamilton cycle.
Let $Mb^{\imath} a^{\jmath}\in N(M)$ for some $\imath\in\ZZ_{10}$ and $\jmath\in\ZZ_p$, where $(\imath,\jmath)\neq(0,0)$.
Suppose that both $\imath\neq0$ and $\jmath\neq0$.
Since $X_{\B}$ is a Cayley graph of a group isomorphic to $\ZZ_{10}$,
it is easy to check that there exists a Hamilton cycle of $X_\B$ containing the edge $\{B_0,B_{\imath}\}$.
Then we get that $d(B_0,B_{\imath})\geqslant2$ as $Mb^{\imath}a^{\jmath}\neq Mb^{\imath}a^{\jmath t^{10}}=Mb^{\imath}a^{\jmath}b^{10}\in \N(M)$, and further, by Proposition~\ref{pro:4}, one yields that $X$ contains a Hamilton cycle, as desired.
Thus, in the following, we may assume that either $\imath=0$ or $\jmath=0$.
Since $X_\B$ contains a Hamilton cycle, we can get a cycle
$${\small\mathcal{A}_1=M\sim Mb^{i_1}\sim Mb^{i_2}\sim\cdots \sim Mb^{i_9}\sim M,}$$
where $\{0,i_1,i_2,\cdots,i_9\}=\{0,1,2,\cdots,9\}$.
According to the connectivity of $X$, there exists an element $Ha^\ell\in \N(M)$ for some $\ell\in\ZZ_p^*$.
Then we get a Hamilton cycle of $X\lg B_0\rg$: 
$${\small\mathcal{A}_2=M\sim Ma^\ell \sim Ma^{2\ell}\sim\cdots \sim Ma^{(p-1)\ell}M\sim M,}$$
which implies $d(B_0)\ge2$.
If $d(B_0)\geq 3$, then by Proposition~\ref{deg3}, $X$ contains a Hamilton cycle, as desired.
So in what follows, we only deal with the case where $d(B_0)= 2$, i.e., $\N(M)\cap B_0=\{Ma^\ell, Ma^{-\ell}\}$.
Since $Ma^{\ell t^{10}}=Ma^\ell b^{10}\in N(M)\cap B_0$, it follows that either $Ma^{\ell t^{10}}=Ma^{\ell}$ or $Ma^{\ell t^{10}}=Ma^{-\ell}$.
Consider the case $Ma^{\ell t^{10}}=Ma^{\ell}$. This implies $a^{\ell(t^{10}-1)}\in M$. Then $t^{10}-1\equiv 0\pmod p$, as $M\cap\lg a\rg=1$, a contradiction.
Consider the case $Ma^{\ell t^{10}}=Ma^{-\ell}$. We get $a^{\ell(t^{10}+1)}\in M$, which implies $t^{10}+1\equiv 0\pmod p$ and so $b^{20}\in \C_G(\lg a\rg)$.
Then we get that $G=\lg a,b\di a^p=b^{20}=1,a^b=a^t,t^{10}\equiv-1\pmod p\rg$ and
$\{Mb^xa^{\ell t^x+y}, Mb^xa^{-\ell t^x+y}\}\subseteq \N(Mb^xa^y)$, where $x\in \ZZ_{10}$ and $y\in\ZZ_p$.
Given the condition $t^{10}\equiv-1\pmod p$, it follows that either $t^4\not\equiv1\pmod p$ or $t^4\equiv1\pmod p$.

Assume $t^4\equiv1\pmod p$.
Then $t^2\equiv-1\pmod p$ and  $\lg a,b^4\rg$ is abelian, which implies $\lg a,b^4\rg=\lg ab^4\rg\cong\ZZ_{5p}$ and $\lg b^4\rg\lhd G$.
Set $\mathcal{C}=\{C_{i',j'} \di C_{i',j'}=M b^{i'}a^{j'}\lg b^4\rg, i'\in\ZZ_2, j'\in\ZZ_p\}$, the set of all orbits of $\lg b^4\rg$.
Consider the quotient graph $X_{\mathcal{C}}$ of order $2p$.
Then $\N(C_{0,0})=\{C_{1,0},C_{0,\ell},C_{0,-\ell}\}$, as the existence of $\mathcal{A}_1$ and $\mathcal{A}_2$.
Since the subgraph of $X_{\mathcal{C}}$ whose edge set is $\{\{C_{0,0},C_{0,\ell}\},\{C_{0,0},C_{0,-\ell}\}\}^G$, is not connected and $X_{\mathcal{C}}$ contains a Hamilton cycle, we get that there exists a Hamilton cycle of $X_{\mathcal{C}}$ containing the edge $\{C_{0,0},C_{1,0}\}$.
If $d(C_{0,0},C_{1,0})\ge2$, then by Proposition~\ref{pro:4}, $X$ contains a Hamilton cycle.
So in what follows, we assume that $d(C_{0,0},C_{1,0})=1$.
Since $X\lg C_{0,0}\cup C_{1,0}\rg=X\lg M\lg b\rg\rg$ contains a Hamilton cycle $\mathcal{A}_1$,
we get that $d(C_{0,0})\ge2$ and $C_{0,0}$ is connected.
If $d(C_{0,0})\ge3$, then $X$ contains a Hamilton cycle by Proposition~\ref{deg3}.
If $d(C_{0,0})=2$, then reset $Y_{i'}= M b^{i'}\lg ab^4\rg$, where $i'=0,1$.
Then $d( Y_0)\geq4$ and $X\lg Y_0\rg$ is connected.
By Proposition \ref{deg3} again, $X$ admits a Hamilton cycle, as desired.

Assume $t^4\not\equiv1\pmod p$.
Since $t^{20}\equiv1\pmod p$, we get
$t^{20}-1=(t^4-1)(t^{16}+t^{12}+t^{8}+t^{6}+t^{4}+1)\equiv0\pmod p$.
By $t^{10}\equiv-1\pmod p$ and $t^4\not\equiv1\pmod p$, we get
$t^{16}+t^{12}+t^{8}+t^{6}+t^{4}+1\equiv 1-t^2+t^4-t^6+t^8\equiv0\pmod p$.
Thus,
\begin{eqnarray}\label{Eq1}
1+t-t^2-t^3+t^4+t^5-t^6-t^7+t^8+t^9=(t+1)(1-t^2+t^4-t^6+t^8)\equiv0~(\bmod ~p).
\end{eqnarray}
By $\mathcal{A}_1$, $\mathcal{A}_2$ and $\{Mb^xa^{\ell t^x+y}, Mb^xa^{-\ell t^x+y}\}\subseteq \N(Mb^xa^y)$,
we get the walks
$${\small\begin{array}{c}
    W_{i_0,j_0}:=Mb^{i_0}a^{j_0}Mb^{i_0}a^{j_0+\ell t^{i_0}}Mb^{i_0}a^{j_0+2\ell t^{i_0}}\cdots Mb^{i_0}a^{j_0+(p-1)\ell t^{i_0}} ~~\text{and}~~\\
    W_{i_0,j_0}^*:=Mb^{i_0}a^{j_0}Mb^{i_0}a^{j_0+(p-1)\ell t^{i_0}}Mb^{i_0}a^{j_0+(p-2)\ell t^{i_0}}\cdots Mb^{i_0}a^{j_0+\ell t^{i_0}},
  \end{array}}$$
where $i_0\in\ZZ_{10}$ and $j_0\in\ZZ_p$.
Thus, in any block $B_{i_0}$, there exists a Hamilton path from $Mb^{i_0}a^{j_0}$ to $Mb^{i_0}a^{j_0+ \delta \ell t^{i_0}}$
and we denote this path by $Mb^{i_0}a^{j_0}\to Mb^{i_0}a^{j_0+ \delta \ell t^{i_0}}$, where $j_0\in\ZZ_p$ and $\delta=\pm1$.
Then we find a Hamilton path $\mathcal{W}$ in $X$:
$${\small\begin{array}{lll}
\mathcal{W}&:=&M\to Ma^\ell\sim Mb^{i_1}a^\ell\to Mb^{i_1}a^{\ell+\delta_{i_1}\ell t^{i_1}} \sim Mb^{i_2}a^{\ell+\delta_{i_1}\ell t^{i_1}}\to Mb^{i_2}a^{\ell+\delta_{i_1}\ell t^{i_1}+\delta_{i_2}\ell t^{i_2}}\sim\\
&&Mb^{i_3}a^{\ell+\delta_{i_1}\ell t^{i_1}+\delta_{i_2}\ell t^{i_2}}\to Mb^{i_3}a^{\ell+\delta_{i_1}\ell t^{i_1}+\delta_{i_2}\ell t^{i_3}+\delta_{i_3}\ell t^{i_3}}\cdots \to Mb^{i_9}a^{\ell+\sum_{s=1}^9\delta_{i_s}\ell t^{i_s}},
\end{array}}$$
where $\delta_1\cdots \delta_9\in\{1,-1\}$.
According to $\mathcal{A}_1$, we get $\ell+\sum_{s=1}^9\delta_{i_s}\ell t^{i_s}=\ell (1+\sum_{s=1}^9\delta_s\ell t^s)$.
By Eq(\ref{Eq1}), we can choose $\delta_1=\delta_4=\delta_5=\delta_8=\delta_9=1$ and $\delta_2=\delta_3=\delta_6=\delta_7=-1$. Then $\ell(1+\sum_{s=1}^9\delta_s\ell t^s)\equiv0\pmod p$,
which implies $Mb^{i_9}a^{\ell+\sum_{s=1}^9\delta_{i_s}\ell t^{i_s}}=Mb^{i_9}\in \N(M)$.
Thus, $\mathcal{W}\sim M$ is a Hamilton cycle, as desired.
\qed

\begin{lemma}\label{Z52}
Using the notations of Theorem \ref{p-blocks}, if $G/K\cong \ZZ_5^2\rtimes\ZZ_{2^\kappa}$ or $\ZZ_2^4\rtimes\ZZ_5$, where $\kappa$ is a  positive integer, then $X$ contains a Hamilton cycle.
\end{lemma}
\demo
Suppose that $G/K\cong \ZZ_5^2\rtimes\ZZ_{2^\kappa}$ with $\kappa$ a positive integer. Note that $\ZZ_p\cong\lg a\rg\unlhd G$ and $\ZZ_p:\ZZ_l\cong\lg a\rg:\lg\sigma\rg=K\leq\AGL(1, p)$. Then $G=(\lg a\rg:\lg\sigma\rg).\ZZ_5^2\rtimes\ZZ_{2^\kappa}\cong (\ZZ_p\rtimes\ZZ_l).\ZZ_5^2\rtimes\ZZ_{2^\kappa}$ and $\C_G(\lg a\rg)\cap K=\lg a\rg$.
Since $G/\C_G(\lg a\rg)$ is cyclic (as $G/\C_G(\lg a\rg)\lesssim \Aut(\lg a \rg)$) and $\ZZ_l\cong\lg\sigma\rg\nleqslant\C_G(\lg a\rg)$, one yields that $\lg a\rg\times N_1\leq\C_G(\lg a\rg)$ for $N_1\cong\ZZ_5^2$.
Further, since $|G|=5^22^\kappa pl$ for $p\geqslant7$, one deduces that $N_1\char C_G(\lg a\rg)$ and $N_1\lhd G$. Then there are $2p$ $N_1$-blocks of size $5$. Of course, $N_1$ acts unfaithfully on some $N_1$-blocks. Applying the Lemma~\ref{unfaithful1}, one yields that $X$ contains a Hamilton cycle, as desired.

Now assume that $G\cong (\ZZ_p\rtimes\ZZ_l).\ZZ_2^4\rtimes\ZZ_5$. By a method completely similar to the above proof, we can deduce that there exits a normal subgroup $N_1\cong \ZZ_2^4$ of $G$ acting unfaithfully on some $N_1$-blocks.
By Lemma~\ref{unfaithful1} again, one yields that $X$ admits a Hamilton cycle, as desired. \qed

\subsection{Cases $r\in \{5, 2\}$}
Throughout this section, we work under  Hypothesis \ref{hypothesis1}. We will specifically deal with the cases where $r=5$ or $2$, corresponding to $m=2p$ or $5p$, respectively. Derived from Lemma \ref{unfaithful1} that if $K$ acts unfaithfully on any $N$-block, then the graph $X$ necessarily contains a Hamilton cycle. Consequently, throughout this section, we will assume that $K$ acts faithfully on every $N$-block.

We begin by describing the structure of the group $G$ in the following.
\begin{lemma}\label{minimal} Suppose $G$ is minimal
transitive group on $\Omega$, of degree $pq$ $($i.e., $|\Omega|=pq$$)$ with $p$-blocks induced by
some normal subgroup of $G$, where $p$ and $q$ are distinct primes. If $G$ is
imprimitive, then $G=P\rtimes Q$, where
 $P\cong \ZZ_p^l$ and $Q\cong \ZZ_{q^k}$ for some positive integers $l$ and $k$,
 and $Q^q\leq K_p$, the kernel of $G$ on $p$-blocks and  $P$ can be
 viewed a $Q/Q^q$-irreducible vector space.
\end{lemma}
\pf Assume that $\Omega$ is a set of size $pq$ with $p$ and $q$ are two distinct primes and $G$ is imprimitive minimal transitive on $\Omega$ with $p$-blocks induced by
some normal subgroup of $G$. Then $G\leq S_p\wr S_q$. The minimal transitivity of $G$ shows that $G\leq\S_p\wr\ZZ_q$. Let $P\in\Syl_p(K_p)$, where $K_p$ is the
kernel of $G$ acting on $p$-blocks. We derive from Proposition \ref{the:main5} that
$G=\N_G(P)K_p$, and so $\N_G(P)$ is transitive on $\O$. By the
minimal transitivity of $G$, one yields that $G=\N_G(P)$, which means $P\unlhd G$. Let
$\sigma \in G$, which permutes cyclicly $p$-blocks.
Then $\sigma^q\leq K_p$. Let $Q\leq \Syl_q(\langle\sigma\rangle)$.  Then $Q^q\leq K_p$ and $P\rtimes Q$ acts transitively on $\O$ with $P\cong \ZZ_p^l$ and $Q\cong \ZZ_{q^k}$. Once again utilizing the minimal transitivity of G, we may conclude that $G=P\rtimes Q$ and $P$ can be
 viewed a $Q/Q^q$-irreducible vector space.  \qed

\begin{lemma} \label{r-rs}
Let $\{r, s, t\}=\{p, 5, 2\}$.
Suppose that $G$ is minimal transitive and imprimitive on $V(X)$. If there exist two normal subgroups $N_1$ and $N_2$ of $G$
 inducing $r$-blocks and $rs$-blocks, respectively, where $N_1\leq N_2$, then $X$ contains
 a Hamilton cycle.
\end{lemma}
\pf Suppose that $G$ is minimal transitive and imprimitive on $V(X)$. Assume that $\{r, s, t\}=\{p, 5, 2\}$, $N_1\unlhd G$ and $N_2\unlhd G$ inducing $r$-blocks and $rs$-blocks respectively. If $r=p$, then by Theorem~\ref{p-blocks}, the graph $X$ contains a Hamilton cycle, as desired.
So in what follows, we assume that $r\neq p$.

(1) Suppose $r=5$. Denote $K_1$ and $K_2$ as the kernels of $G$ acting on its systems of $5$-blocks and $5s$-blocks, respectively.
If $K_1$ is unfaithful on some $5$-block, then by Lemma~\ref{unfaithful1}, the graph $X$ must contain a Hamilton cycle. Therefore, it is sufficient to consider the case where $K_1$ acts faithfully on every $5$-block.

Suppose $s=p$ (so $t=2$). Then $G/K_1\leq \Aut(X_{K_1})$ and $|V(X_{K_1})|=2p$. Let $P\in\Syl_p(G)$. By Lemma~\ref{minimal}, one yields that $G/K_1=PK_1/K_1\rtimes \langle \sigma \rangle K_1/K_1$, where
$\langle \sigma \rangle$ is a $2$-group. Noting that $K_1\leq\S_5$,
$P\in\Syl_p(G)$ for $p\geqslant 7$ and $PK_1=K_1\rtimes P$, one yields that $\Aut(K_1)$ is a $\{2, 3, 5\}$-group (checked by Magma), and so $PK_1=P\times K_1$.
Consequently, ${\rm{Syl}}_p(PK_1)=\{ P\}$. Followed
from $PK_1\lhd G$ that $P\lhd G$, and so $P$ induces $p$-blocks of
$G$.
By Theorem~\ref{p-blocks}, the graph $X$ contains a Hamilton cycle, as desired.

Suppose $s=2$ (so $t=p$). Building upon Lemma~\ref{minimal} and the preceding argument, we obtain $G/K_1\cong \ZZ_2^l\rtimes \ZZ_p$ and thus $|G|=|K_1|2^lp$ for some positive integer $l$.
Then $G$ can be expressed as $K_1.L/K_1\rtimes PK_1/K_1$ where $L/K_1\cong\ZZ_2^l$ and $P\cong\ZZ_p$ is
a Sylow $p$-subgroup of $G$.
Note that $p\geqslant7$ and $K_1\leq\S_5$. Let $Q\in \Syl_5(G)$.
By Proposition \ref{the:main5}, we have $G=\N_G(Q)K_1$.
It follows that $\N_G(Q)$ is transitive on $\Omega$ and thus $G=\N_G(Q)$, as $G$ is a minimal transitive subgroup on $V(X)$ and $\N_G(Q)\leq G$. Consequently, $Q\unlhd G$ and $Q\unlhd K_1$, implying that $\Syl_5(K_1)=\Syl_5(G)=\{ Q\}$.
Now, consider the conjugation action of $L/K_1$ on the Sylow $5$-subgroups of $K_1$.
Since  $Q^{(L/K_1)}=Q$, we deduce that $Q\cong\ZZ_5$, $\C_{K_1}(Q)=Q$ and $QL/K_1\leq \C_G(Q)$ (by the Frattini-Argument).
Therefore, $\C_G(Q)/Q\cong \C_G(Q)K_1/K_1\geqslant L/K_1$, which implies that either $\C_G(Q)=Q\times L_1$ or $\C_G(Q)=Q\times (L_1\rtimes P)$, where $\ZZ_2^l\cong L_1\leq\C_G(Q)$.
In both cases, $L_1\lhd G$, and thus $L_1$ induces $2$-blocks as $|V(X)|=10p$.
We now focus on the case where $L_1$ acts faithfully on each $2$-block, which implies $l=1$.
Consequently, $G/K_1=L/K_1\times PK_1/K_1$, as $\ZZ_2.\ZZ_p=\ZZ_2\times\ZZ_p$.
From $K_1\leq\S_5$, we deduce that $PK_1=P\times K_1$, and thus $P\char PK_1$. Consequently, $P\lhd G$, and the normal subgroup $P$ induces a system of $p$-blocks on $V(X)$. By Theorem~\ref{p-blocks}, the graph $X$ must contain a Hamilton cycle, as desired.

(2) Suppose $r=2$.
Let $K_1$  and $K_2$ denote the kernels of the action of $G$ on its systems of $2$-blocks and $2s$-blocks, respectively. Following the previous approach, we may restrict our consideration to the case where $K_1\cong \ZZ_2$ is faithful on each $2$-block.

Suppose $s=p$ (so $t=5$).  Then $G/K_1$ is a transitive permutation group of degree $5p$.
Let $P\in \Syl_p(G)$. By Lemma~\ref{minimal}, we may set $G/K_1=PK_1/K_1\rtimes \langle \sigma \rangle
K_1/K_1$, where $\langle \sigma \rangle$ is a 5-group.
Since $K_1 \cong \ZZ_2$, it follows that $PK_1 = P \times \ZZ_2$, and $P \lhd G$. Consequently, $P$ induces $p$-blocks of $G$.
Then, by Theorem~\ref{p-blocks}, $X$ contains a Hamilton cycle, as desired.

Suppose $s=5$ (so $t=p$).
Then reset $G/K_1=QK_1/K_1\rtimes PK_1/K_1$, where $Q\in \Syl_5(G)$ and $Q\cong \ZZ_5^l$ for some positive integer $l$ (by Lemma \ref{minimal}).
Here, $Q \lhd G$ and thus $Q$ induces $5$-blocks. If $Q$ is unfaithful on a $5$-block, then $X$ contains a Hamilton cycle by Lemma \ref{unfaithful1}. Now, assume that $Q$ acts faithfully on each $5$-block. Then $Q\lesssim \S_5$, imiplying $Q\cong \ZZ_5$.
Consequently, $P \lhd G$, and it follows that there exists a normal subgroup in $G$. By Theorem~\ref{p-blocks} again, $X$ contains a Hamilton cycle, as desired.
\qed

\begin{theorem} \label{q=5}
$X$ is  Hamiltonian if $r=5$.
\end{theorem}
\demo
Now $|B_i|=r=5$ for $i=0, 2,\dots, m-1$ and $m=2p$. Since $K$ acts faithfully on each block $B_i$ and $G$ is a minimal transitive group on $V(X)$, we obtain that $\ZZ_5\lesssim K\lesssim\S_5$ and $G/K$ is a minimal transitive permutation group of degree $2p$. Let $L\in\Syl_5(K)$. Then $L\cong\ZZ_5$ is transitive on $B_i$ for $0\leq i\leq 2p-1$, and by Proposition \ref{the:main5}, $G=\N_G(L)K$ and $\N_G(L)\leq G$ is transitive on $V(X)$. This implies $G=\N_G(L)$ by the minimal transitivity of $G$, and hence $L\unlhd G$ and $L\unlhd K$, meaning $\Syl_5(K)=\{L\}$.
Then $(\C_G(L)K)/K\trianglelefteq G/K$.

By Lemma \ref{r-rs}, we only need to consider the case where $G/K$ is quasprimitive on $V(X_{\mathcal{B}})$ of degree $2p$, where $\mathcal{B}$ is a complete $K$-block system of $G$. By the minimal transitivity of $G$, we assume that $T:=G/K=\soc(G/N)$.
By Propositions \ref{2p} and \ref{quasiprimitive}, we deduce that $T$ is a nonabelian simple group.
Then either $\C_G(L)K=K$ or $\C_G(L)K=G$, as $(\C_G(L)K)/K\trianglelefteq G/K$. If $\C_G(L)K=K$, then $\C_G(L)\leq K$.
It follows that $\C_G(L)=L$, as $K\lesssim \S_5$ and $\C_{S_5}(\ZZ_5)=\ZZ_5$. This contradicts the fact that $G/\C_G(L)\lesssim \Aut(\ZZ_5)\cong \ZZ_4$ and $10p\mid |G|$. Thus, we conclude that $\C_G(L)K=G$. Then either $\C_G(L)$ induced $5p$-blocks on $V(X)$ or $\C_G(L)$ is transitive on $V(X)$ and hence $G=\C_G(L)$. If $\C_G(L)$ induced $5p$-blocks on $V(X)$, then $X$ admits a Hamilton cycle by Lemma \ref{r-rs}, as $L$ induces $5$-blocks on $V(X)$.

Now, set $\C_G(L)=G$. Then $L=K$, as $K\leq \C_G(L)$ and $K\lesssim S_5$.
Moreover, $\C_G(L)/L=G/K=T$ is quasiprimitive on $V(X_{\mathcal{B}})$ of degree $2p$.
If $T$ is primitive, then by Proposition \ref{direct product}, the graph $X$ contains a Hamilton cycle.
Note that $T$ is quasiprimitive but not primitive. Thus, by Lemma \ref{quasiprimitive}, either $T\cong\M_{11}$ with $p=11$ or $\PSL(n, r^f)$ for primes $n$ and $r$, a positive integer $f$, with $(n, r^f-1)=1$ and $p=\frac{r^{nf}-1}{r^f-1}$.
After carefully examining the Schur multiplier of $T$, we conclude that $5\nmid Mult(T)$.
Consequently, $G=L\times T\cong\ZZ_5\times T$, where either $T\cong \M_{11}$ with $p=11$ or $T\cong \PSL(d, e)$ with $p=\frac{e^{d}-1}{e-1}$, where $d$ and $e$ are two odd primes satisfying $(d,e-1)=1$.
If $T\cong \M_{11}$, then a Hamilton cycle of $X$ can be obtained using Magma.
If $T\cong \PSL(d,e)$, then a Hamilton cycle of $X$ follows from Lemma \ref{quasi-pri-2p}.
\qed

\begin{theorem} \label{q=2}
$X$ is  Hamiltonian if $r=2$.
\end{theorem}
\demo
Suppose that $r=|B_i|=2$ for $i=1, 2, \dots, m-1$, i.e., $m=5p$. Since $K$ is faithful on $B_i$ for $0\leq i\leq 5p-1$, it follows that $K\cong\ZZ_2$ and $G/K$ is a minimal transitive permutation group of degree $5p$.
By Lemma \ref{r-rs}, we need only consider the case where  $G/K$ is quasiprimitive on $V(X_K)$, i.e., $G/K$ is a quasiprimitive permutation group of degree $5p$.
Then, Propositions \ref{5p} and \ref{quasiprimitive} provide all possible candidates for $G/K$, implying that $\soc(G/K)$ is an almost simple group. Applying the minimal transitivity of $G$, we deduce that $G/K=\soc(G/K)$ is a nonabelian simple group, denoted by $T$.
If $G/K$ is $2$-transitive on $V(X_K)$, then $X_K\cong\K_{5p}$ and the valency of $X_K$ is $5p-1$, which exceeds  $\frac{10p}{3}$. Consequently, the valency of $X$ also exceeds $\frac{10p}{3}$.
By Proposition \ref{pro:2}, the graph $X$ contains a Hamilton cycle.
Now, assume that $G/K$ is not 2-transitive on $V(X_K)$.

Assume that $G/K$ is primitive on $V(X_K)$. By Proposition \ref{5p}, if $G/K\cong\M_{11}$, $\PSL(2, 11)$, $\A_{11}$ or $\PSL(5, 2)$, then $T$ contains a transitive metacyclic group $\ZZ_p\rtimes\ZZ_5$. Using the minimal transitivity of $G/K$, one deduces that $T=G/K=\ZZ_p\rtimes\ZZ_5$. This contradicts the simplicity of $T$. We now analyze the remaining cases from Proposition \ref{5p}. Let $\a\in B_0$.
Then $G_{B_0}=K\rtimes G_\a$, meaning that $K$ has a complement in $G_{B_0}$.
Combing these facts: $K\leq G_{B_0}\leq G$, $|G:G_{B_0}|=|V(X_K)|=5p$ for a prime $p \geq 7$, and $K\cong\ZZ_2$, we deduce that $(|K|, |G:G_{B_0}|)=(2, 5p)=1$. By Proposition \ref{complement}, $K$ has a complement in $G$, denoted by $\overline{K}$.
Since $\overline{K} \cong G/K = T$, it follows that $G \cong K \times T$.
Since $T = G/K$ is primitive on $V(X_K)$, it acts primitively on the set $[T : T_{B_0}]$. Therefore, by Proposition \ref{direct product}, the graph $X$ contains a Hamilton cycle, as desired.

Suppose that $G/K$ is quasiprimitive but not primitive on $V(X_K)$. Proposition \ref{quasiprimitive} shows that $G/K\cong\SL(2,2^{2^s})$ for $p=2^{2^s}+1$ and $5|2^{2^s}-1$, $\PSL(d,e)$ for some odd primes $d$ and $e$, $p=\frac{e^d-1}{e-1}$ and $5\mid e-1$, or $\PSL(2, 11)$ for $p=11$.
The latter former coincide with the corresponding situation above. For the first form, since the Schur multiplier of $G/K=T\cong\SL(2, 2^{2^s})$ is $1$ and $N\cong\ZZ_2$, we have $G=K.T=K\times T\cong\ZZ_2\times\SL(2, 2^{2^s})$.
This implies that $T$ has $2$ blocks in $V(X)$. Observe that $T$ is faithful and quasiprimitive on both of its blocks. By Proposition \ref{5p},  $T$ is not primitive on every $T$-block. By Proposition \ref{quasi-pri-2}, the graph $X$ contains a Hamilton cycle, as desired.
If $G/K\cong\PSL(d, e)$ for some odd primes $d$ and $e$ with $(d,e-1)=1$, $p=\frac{e^d-1}{e-1}$ and $5\mid e-1$, then  $G=K.T=K\times T\cong\ZZ_2\times\PSL(d,e)$, and $T$ has $2$ blocks in $V(X)$. It is clear that $T$ is faithful and quasiprimitive on both $T$-blocks. By Proposition \ref{5p} again, $T$ is not primitive on any $T$-orbit.
By Lemma \ref{quasi-pri-1}, the graph $X$ contains a Hamilton cycle, as desired.
\qed

\section{Cases $r\in\{10,5p,2p\}$ }\label{10-blocks}
In this section, we will establish three theorems that handle the cases where $K$ induces blocks of length $r$ (with $r \in \{10, 5p, 2p\}$), under Hypothesis \ref{hypothesis1}.

\begin{theorem}\label{l=10}
$X$ is  Hamiltonian if $r=10$.
\end{theorem}
\demo
By Hypothesis \ref{hypothesis1}, $\mathcal{B}=\{B_0,B_1,B_2,\dots,B_{p-1}\}$ is the set of all $K$-orbits, where $K \trianglelefteq G$ and $G/K \cong \ZZ_p$.
In particular, $G$ is a minimal transitive subgroup on $V(X)$, and $G/K$ is a minimal transitive subgroup on $V(X_{\mathcal{B}})$.
Consequently, $G=K\lg \tau \rg$ for some element $\tau$ of order $p^e$, and the subgroup $\lg \tau \rg$ acts transitively on $\mathcal{B}$. This implies that the quotient graph $X_{\mathcal{B}}$ contains a cycle of length $p$. Without loss of generality, we may assume that $B_{i+1} \in N(B_i)$ for $i \in \mathbb{Z}_p$.

Let $N_0$ be a minimal normal subgroup of $G$ contained in $K$. If $N_0$ induces $N_0$-blocks of size $5$ or $2$, then by Theorems \ref{q=5} and \ref{q=2}, $X$ admits a Hamilton cycle, as desired. Therefore, for the remainder of the argument, we assume that every minimal normal subgroup of $G$ contained in $K$ induces blocks of size 10. Without loss of generality, let $N_0$ be one such minimal normal subgroup of $G$. Then the group $N_0\lg \tau \rg$ is transitive on $V(X)$, which implies that $G = N_0\lg \tau \rg$.

If $N_0$ is unfaithful on some $K$-orbit, then, by Lemma~\ref{unfaithful2}, $X$ admits a Hamilton cycle. Therefore, for the remainder of the argument, we assume that $N_0$ is faithful on each $K$-orbit. Then, by Proposition \ref{squarefree}, $N_0$ embeds into $S_{10}$ (i.e., $N_0 \lesssim S_{10}$) and is simple and quasiprimitive on all $K$-orbits.

Note that $G = N_0\lg\tau \rg$, where $\lg\tau \rg\cong \ZZ_{p^e}$, $N_0 \lesssim S_{10}$, and $p \geq 7$ is prime. If $(|N_0|, p) = 1$, then $(|\Aut(N_0)|, p)=1$ (verified using Magma), and the `NC'-theorem implies that $G \cong N_0 \times \langle \tau \rangle$. Consequently, there exists a normal subgroup of $G$ that induces $p$-blocks in $V(X)$.
By Theorem~\ref{p-blocks}, the graph $X$ contains a Hamilton cycle, as desired.
Now suppose that $(|N_0|, p) \neq 1$ for a prime $p \geq 7$ and $N_0 \lesssim S_{10}$. This implies $(|N_0|, p) = 7$, meaning that $7\mid|N_0|$ and $p = 7$. Thus, $\ZZ_7\lesssim N_0$ is a permutation group of degree $10$. Consulting \cite[p. 48]{Atlas}, we conclude that $N_0 \cong A_{10}$, and $\Aut(N_0) \cong S_{10}$, since $N_0$ is simple and possesses a subgroup of index $10$. With $p = 7$, we have $G \lesssim S_{70}$. This implies $\tau\in\S_{70}$, so its order is either $7$ or $7^2$. Using Magma, we confirm that $X$ contains a Hamilton cycle, as desired.\qed

\vskip 5mm

\begin{theorem}\label{l=5p}
$X$ is  Hamiltonian if $r=5p$.
\end{theorem}
\demo Now, let $X$ be a graph of order $10p$, where $p \geqslant 7$ is a prime, and $G$ is a minimal transitive subgroup of $\Aut(X)$. Moreover, let $N \unlhd G$ such that all the orbits of $N$ on $V(X)$ form a complete block system of $G$, denoted by $\B=\{B_0, B_1\}$, where $|B_i| = 5p$ for $i = 0, 1$. 
In particular, $G/K$ is a minimal transitive subgroup of $\Aut(X_\B)$ acting on $V(X_\B)$, where $K$ is the kernel of this action. Then there exists a positive integer $e$ and an element $\tau \in G$ such that $\o(\tau) = 2^e$, $G = K \langle \tau \rangle$, and $\tau$ swaps the two blocks $B_0$ and $B_1$. Note that $X[B_0, B_1]$ is a generating subgraph of $X$.

Suppose that $K$ is unfaithful on an $N$-block, say $B_0$. Let $T_0 = K_{(B_0)}$ be the kernel of $K$ acting on the set $B_0$, and let $T=T_0\times T_0^\tau$.
Then $T \unlhd G$ and $T\leq K$. If $T$ is transitive on $B_0$ or $B_1$, then $X[B_0, B_1] \cong \K_{5p, 5p}$, and so it contains a Hamilton cycle by Proposition \ref{pro:2}.
Therefore, Theorem \ref{l=5p} holds, as $X[B_0, B_1]$ is a generating subgraph of $X$.
Now assume that $T$ is not transitive on $B_0$ and $B_1$, and it induces blocks of length $r'$. Then $r' \in\{5, p\}$.
We deduce from Theorems \ref{p-blocks} and \ref{q=5} that $X$ admits a Hamilton cycle, as desired.
So, in what follows, we assume that $K$ is faithful on all $N$-blocks.

Assume, without loss of generality, that $K$ is primitive on an $N$-block, say $B_0$. Then $K$ is a primitive permutation group of degree $5p$ for $p \geqslant 7$. Let $K' = \soc(K)$. By Proposition \ref{5p}, we have that $(K', p)=(\A_{5p}, p)$, $(\PSL(d,q), \frac{q^d-1}{5(q-1)})$, $(\P\Omega(5,4), 17),(\PSU(3,4), 13)$,
$(\A_{11}, 11),(\A_7, 7),(\PSL(4,2), 7),(\PSL(5,2), 31),(\PSp(4,4), 17),(\P\Omega^+(6,2),7),(\PSL(2, 5^2)$,
$13)$, $(\M_{11}, 11)$ or $(\PSL(2, 11), 11)$.
Moreover, $K$ is $2$-transitive in the first four cases, and $K'$ has a metacyclic transitive subgroup isomorphic to $\ZZ_p : \ZZ_5$ in the following cases: $K' \in \{\A_{11}, \PSL(5, 2), \M_{11}, \PSL(2, 11)\}$. Take $u_0 \in B_0$.
Since the stabilizer $K_{u_0}$ is maximal in $K$ if and only if $K_{u_0^\tau}$ is maximal in $K$, we deduce that $K$ is primitive on $B_1$. Then $K' \char K \unlhd G$ and $K' \unlhd G$. It follows from $K$ acting primitively on $B_0$ and $B_1$ that $K'$ is transitive on both $B_0$ and $B_1$. Set $G_0 = K' \lg \tau \rg$ and $K_0 = K'\lg \tau^2 \rg$. Then $G_0$ is transitive on $V(X)$, $K_0\unlhd G_0$, and $G_0/K_0 \cong \ZZ_2$. Further, the minimal transitivity of $G$ implies that $G = G_0$.

\vskip 3mm
{\it Claim 1:
 If $K'$ has a metacyclic transitive subgroup isomorphic to $\ZZ_p:\ZZ_5$, then the graph $X$ admits a Hamilton cycle.}
\vskip 3mm
Further, since $\tau$ swaps $B_0$ and $B_1$, we have $P^\tau = P^g$ for some $g \in K'$. This implies $\tau g^{-1} \in \N_G(P)$, and consequently, $\N_G(P) \geq \lg L, \tau g^{-1} \rg$. In particular, $\N_G(P)$ is transitive on $V(X)$.
The normal Sylow $p$-subgroup $P$ of $\N_G(P)$ induces $p$-blocks on $V(X)$. By the minimal transitivity of $G$ on $V(X)$, we conclude that $G = \N_G(P)$. By Theorem~\ref{p-blocks}, we deduce that $X$ admits a Hamilton cycle.

\vskip 3mm
{\it Claim 2:
If $K'$ has no transitive metacyclic subgroup on $B_i$ for $i=0,1$, then the graph $X$ admits a Hamilton cycle.}
\vskip 3mm
If $(K',p)\in\{(A_{5p},p), (\PSL(d,q),\frac{q^d-1}{5(q-1)}), (\P\Omega(5,4),17), (\PSU(3,4), 13)\}$, that is $K$ is $2$-transitive on both $B_0$ and $B_1$.
Consequently, the bipartite graph $X[B_0, B_1]$ is isomorphic to $K_{5p,5p}$, $K_{5p,5p}$-$5pK_{1,1}$, or $5pK_{1,1}$.
If $X[B_0, B_1] \cong K_{5p,5p}$ or $K_{5p,5p} - 5pK_{1,1}$, then the minimum degree $d(X)\geq 5p-1\geq \frac{10p}{3}$. By Proposition \ref{pro:2}, $X$ admits a Hamilton cycle.
If $X[B_0, B_1] \cong 5pK_{1,1}$, then, due to the connectivity of $X$, there must exist an edge within $X[B_0]$. Since $K$ is $2$-transitive on $B_0$, it follows that $X[B_0] \cong K_{5p}$. Therefore, $d(X) \geq 5p - 1 \geq \frac{10p}{3}$, and applying Proposition \ref{pro:2} again, we deduce that $X$ contains a Hamilton cycle.
Now, if  $(K',p)\in \{(A_7,7)$, $(\PSL(4,2), 7)$, $(\PSp(4,4), 17)$, $(\P\Omega^+(6,2),7)$, $(\PSL(2,5^2), 13)\}$, then, using Magma, one can verify that the graph $X$ contains a Hamilton cycle.
Hence, Claim 2 holds.

Suppose that $K$ is imprimitive and quasiprimitive on both $N$-blocks. By the arguments of Claim 1, we deduce that the graph $X$ admits a Hamilton cycle whenever $\soc(K) = \PSL(2,11)$. Combining this with Propositions \ref{quasi-pri-2} and \ref{quasiprimitive} and Lemma \ref{quasi-pri-1}, we conclude that $X$ contains a Hamilton cycle.

Suppose that $K$ is not quasiprimitive on an $N$-block, say $B_0$, without loss of generality. Then there exists a maximal intransitive normal subgroup $K'$ of $K$ acting on $B_0$, and $K'$ induces blocks of length $5$ or $p$ on $B_0$. Note that $G = K \langle \tau \rangle$ and $G/K \cong \ZZ_2$, with $\tau$ swapping $B_0$ and $B_1$. Then $K'K'^\tau \leq K$ and $K'K'^\tau \unlhd G$. Since $K'$ is a maximal intransitive normal subgroup of $K$, we have that either $K' = K'^\tau$ or $K'K'^\tau$ is transitive on $B_0$.
For the former case, i.e., $K' = K'^\tau$, it follows that $K' \trianglelefteq G$. Moreover, $K'$ induces $r'$-blocks for some $r' \in \{5, p\}$ on $B_0$. By Theorems \ref{p-blocks} and \ref{q=5}, the graph $X$ contains a Hamilton cycle.

Suppose that $K'K'^\tau$ is transitive on $B_0$. If $K'$ induces blocks of length $5$ on $B_0$, then $K'K'^\tau/K'$ is a transitive permutation group on $p$ points. Hence, $p \mid |K'K'^\tau/K'|$ and $p \mid |K'|$. However, since $K'$ acts faithfully on $B_0$, we have $K' \lesssim S_5^p$, which contradicts $p \mid |K'|$ (as $p \geq 7$). Therefore, $K'$ must induce $p$-blocks on $B_0$, forming a set $\mathcal{D}$.
Let $P$ be a Sylow $p$-subgroup of $K'$ and $\mathcal{N}_p$ be the set of all $p$-elements of $K$ (elements whose order is a power of $p$).
Since $P^{\tau}\leqslant K'^{\tau}\leqslant K$ and $G=K\lg\tau\rg$, for any $g \in G$, we have $P^g\leqslant \lg P, \mathcal{N}_p \rg$. Thus, $\lg P, \mathcal{N}_p \rg\unlhd G$.
Clearly, $P$ induces the same $p$-blocks as those of $K'$ acting on $B_0$.
On the other hand, since $|\mathcal{D}|=5$ and $p\geq7$, any element in $\mathcal{N}_p$ fixes each block in $\mathcal{D}$.
Consequently, $\lg P, \mathcal{N}_p \rg$ induces $p$-blocks in $V(X)$.
By Theorem~\ref{p-blocks},  $X$ contains a Hamilton cycle, as desired.\qed

\vskip 5mm

\begin{theorem}\label{l=2p}
$X$ is  Hamiltonian if $r=2p$.
\end{theorem}
\demo
By Hypothesis \ref{hypothesis1}, we deduce that $\B =\{B_0,B_1, B_2, B_3, B_4\}$ is the set of all $N$-orbits, where $N \unlhd G$ and $G/K \cong\ZZ_5$. In particular, $G$ and $G/K$ are minimal transitive subgroups acting on $V(X)$ and $V(X_N)$, respectively.
Then, $G=K\lg\t\rg$ for some $\t$ of order $5^{e'}$, where $e'$ is a positive integer, and $\lg\t\rg$ acts transitively on $\B$. Consequently, the quotient graph $X_N$ has a cycle of length $5$. Without loss of generality, we assume that $B_{i+1} \in \N(B_i)$ for $i \in \ZZ_5$.

Let $N_0$ be a minimal normal subgroup of $G$ contained in $K$.
If $N_0$ induces $N_0$-blocks of size $p$ or $2$, then by Theorems \ref{p-blocks} and \ref{q=2}, $X$ contains a Hamilton cycle, as desired. So in what follows, we assume that every minimal normal subgroup of $G$ contained in $K$ induces blocks of size $2p$. Then $N_0\lg \tau\rg$ acts transitively on $V(X)$, implying $G=N_0\lg \tau \rg$.

Suppose that $N_0$ is unfaithful on an $N$-orbit.
By Lemma \ref{unfaithful2}, $X$ admits a Hamilton cycle.
Thus, in the following, we assume that $N_0$ is faithful on all $N$-orbits.
Further, by Proposition \ref{squarefree}, $N_0$ is simple and acts quasiprimitively on all $N$-orbits.

Suppose that $N_0$ is primitive on an $N$-orbit, say $B_0$, without loss of generality. Let $u_0 \in B_0$.
The stabilizer ${N_0}{u_0}$ of $u_0$ in $N_0$ is maximal in $N_0$ if and only if ${N_0}{u_0^{\tau'}}$ is maximal in $N_0$ for any $\tau'\in\lg\tau\rg$.
This implies that $N_0$ is primitive on $B^{\tau'}$ and hence, primitive on all $N$-orbits. Let $K'=\soc(N_0)$.
Then $K'\char N_0\unlhd G$ and $K'\lhd G$.
Combining this with Proposition \ref{2p} and $p \geq 7$, we conclude that
$(K', 2p)\in\{(\M_{22}, 22)$, $(\A_{2p}, 2p)$, $(\PSL(2, r'^{2^\imath}), {r'}^{2^\imath}+1)\}$, where $r'$ is an odd prime and $i$ is a positive integer.
Clearly, $K'$ is 2-transitive on $B_i$, and thus, $X[B_i, B_{i+1}]\cong\K_{2p, 2p}$, $\K_{2p, 2p}-2p\K_{1,1}$ or $2p\K_{1,1}$, where $i\in\ZZ_5$.
For the first two cases, $d(B_i, B_j) \geq 2p - 1$ for $B_j\in X_\B'(B_i)$. Since $d(X_\B')\geq2$ and $p\geq7$, we deduce that $d(X) \geq 2(2p - 1) \geq \frac{10p}{3}$.
By Proposition \ref{pro:2}, $X$ admits a Hamilton cycle.
Now, suppose that $X[B_i, B_{i+1}] \cong 2p K_{1,1}$ for $i \in \ZZ_5$.
Then $X\lg B_i\rg$ is not an empty graph, as $X$ is connected, where $i\in\ZZ_5$.
Hence, $X\lg B_i \rg\cong \K_{2p}$, as $K'$ is $2$-transitive on each $B_i$ for $i \in \ZZ_5$.
It is then straightforward to construct a Hamilton cycle in the graph $X$, as desired.

Suppose that $N_0$ is imprimitive and quasiprimitive on all $N$-orbits. Let $K'=\soc(N_0)$, by the same arguments as in the previous paragraph.
By Proposition \ref{quasiprimitive}, either $K'\cong \M_{11}$ or $K'\cong \PSL(d,r')$, where $(d, r'-1)=1$, $d$ and $r'$ are both odd primes, and $p=(r'\sp d-1)/(r'-1)$.
For the former case, $K' \cong M_{11}$. By the minimal transitivity of $G$, we have $\M_{11}\times\ZZ_5\lesssim G\lesssim \Aut(\M_{11})\times\ZZ_5$.
Then $G$ has a normal subgroup that induces 5-blocks. By Theorem \ref{q=5}, $X$ admits a Hamilton cycle, as desired.
Now, suppose that $K'\cong \PSL(d, r')$, where $(d, r'-1)=1$, $d$ and $r'$ are odd primes, and $p=(r'\sp d-1)/(r'-1)$. Then, by Lemma \ref{quasi-pri-2p},  $X$ contains a Hamilton cycle.
\qed

\vskip 5mm

Finally, we conclude this paper with the proof of Theorem \ref{main}.

\f{\bf Proof of Theorem \ref{main}:}
Let $X$ be a $G$-vertex transitive graph of order $10p$, where $p$ is a prime.
First, consider the case where $p\leq5$. By \cite{DZ2025,KM08,KM09,DM87}, it follows that $X$ is either a truncation of the Petersen graph or a Hamiltonian graph.
Next, assume that $p\geq7$.
If $X$ is a quasiprimitive (resp. primitive) graph, then $X$ is a Hamiltonian graph by \cite[Proposition 4.1]{KMZ12} and \cite[Theorem 1.1]{DLY-25} (resp. \cite[Theorem 1.1]{DTY}).
So in what follows, we focus on the case where $X$ is a non-quasprimitive graph.
In this scenario, there exists a normal subgroup $N$ of $G$ that induces nontrivial blocks for $G$.
Let $r$ denote the size of each such $N$-block. Our analysis shows $r\in\{2,5,p,10,2p,5p\}$.
We now apply the appropriate theorem for each possible value of $r$: Theorems \ref{p-blocks}, \ref{q=5}, \ref{q=2}, \ref{l=10}, \ref{l=5p}, and \ref{l=2p} collectively imply that $X$ is  a Hamiltonian graph in all these cases. This completes the proof.
\qed

{\it Email address:} Huye Chen: chenhy280@gxu.edu.cn

Jingjian Li: lijjhx@gxu.edu.cn

Hao Yu: haoyu@gxu.edu.cn

\end{document}